\documentclass[11pt]{article}
\usepackage{amsfonts}
\usepackage{color}
\usepackage{amsmath,amssymb,comment}
\usepackage[numbers,sort&compress]{natbib}
\usepackage{mathtools}
\newtheorem{theo}{Theorem}[section]
\newtheorem{lem}[theo]{Lemma}

\newtheorem{defi}[theo]{Definition}

\newcommand{\mysection}[1]{\section{#1} \setcounter{equation}{0}}
\newcommand{\proof}{{\sc Proof.} \quad}
\newcommand{\proofc}{{\sc Proof} \ }
\newcommand{\be}{\begin{equation} \label}
	\newcommand{\ee}{\end{equation}}
\newcommand{\bea}{\begin{eqnarray}\label}
	\newcommand{\eea}{\end{eqnarray}}
\newcommand{\bas}{\begin{eqnarray*}}
	\newcommand{\eas}{\end{eqnarray*}}
\newcommand{\bit}{\begin{itemize}}
	\newcommand{\eit}{\end{itemize}}
\newcommand{\qed}{\hfill$\Box$ \vskip.2cm}
\newcommand{\qedd}{\hfill$\Box$ \vskip.2cm}
\newcommand{\nn}{\nonumber}
\newcommand{\R}{\mathbb{R}}
\newcommand{\N}{\mathbb{N}}
\newcommand{\pO}{\partial\Omega}

\newcommand{\eps}{\varepsilon}

\newcommand{\wto}{\rightharpoonup}
\newcommand{\wsto}{\stackrel{\star}{\rightharpoonup}}
\newcommand{\hra}{\hookrightarrow}
\newcommand{\io}{\int_\Omega}
\newcommand{\na}{\nabla}
\newcommand{\Del}{\Delta}
\newcommand{\del}{\delta}
\newcommand{\al}{\alpha}
\newcommand{\lam}{\lambda}

\newcommand{\sig}{\sigma}
\newcommand{\pa}{\partial}
\newcommand{\bom}{\overline{\Omega}}
\newcommand{\Om}{\Omega}

\newcommand{\ov}{\overline}

\newcommand{\wh}{\widehat}

\newcommand{\hs}{\hspace*}

\newcommand{\vp}{\varphi}
\newcommand{\lbal}{\left\{ \begin{array}{l}}
	\newcommand{\lball}{\left\{ \begin{array}{ll}}
		\newcommand{\ear}{\end{array} \right.}

	\newcommand{\cred}{\color{red}}

	\newcommand{\abs}{\\[5pt]}
	\newcommand{\Abs}{\\[5mm]}
	
	\newcommand{\adb}{\allowdisplaybreaks}
	\newcommand{\tm}{T_{max}}
	
	\newcommand{\tmN}{T_{max,N}}

	\newcommand{\ou}{\ov{u}}

	\newcommand{\F}{{\mathcal{F}}}
	
\newcommand{\Ltwo}{L^2(\Omega)}
\newcommand{\Honezero}{W_0^{1,2}(\Omega)}		
\newcommand{\Htwoone}{W^{2,2}(\Omega) \cap W_0^{1,2}(\Omega)}
\newcommand{\change}[1]{\textcolor{black}{#1}}		

\begin{document}
\adb
%
%
\title{$L^\infty$ blow-up in the Jordan--Moore--Gibson--Thompson equation}
\author{
Vanja Nikoli\'c\footnote{vanja.nikolic@ru.nl}\\
{\small Department of Mathematics,
	Radboud University}\\
{\small 6525 AJ Nijmegen, The Netherlands} 
\and
Michael Winkler\footnote{michael.winkler@math.uni-paderborn.de}\\
{\small Universit\"at Paderborn, Institut f\"ur Mathematik}\\
{\small 33098 Paderborn, Germany}
}
\date{}
\maketitle
\begin{abstract}
\noindent 
The Jordan--Moore--Gibson--Thompson equation
\bas
	\tau u_{ttt} + \al u_{tt} = \beta \Del u_t + \gamma \Del u + (f(u))_{tt}
\eas
is considered in a smoothly bounded domain $\Om\subset\R^n$ with $n\le 3$, where $\tau>0,\beta>0,\gamma>0$, and $\al\in\R$.\abs
Firstly, it is seen that under the assumption that $f\in C^3(\R)$ is such that $f(0)=0$, gradient blow-up phenomena
cannot occur in the sense that for any appropriately regular initial data, within a suitable framework of strong solvability, 
an associated Dirichlet type initial-boundary value problem admits a unique solution $u$ on a maximal time interval $(0,\tm)$ which
is such that
\bas
	\mbox{if $\tm<\infty$, \quad then \quad} \limsup_{t\nearrow\tm} \|u(\cdot,t)\|_{L^\infty(\Om)}=\infty.
	\qquad \qquad (\star)
\eas
This is used to, secondly, make sure that if additionally $f$ is convex and grows superlinearly in the sense that
\bas
	f''\ge 0
	\mbox{ on $\R$,} 
	\qquad
	\frac{f(\xi)}{\xi} \to +\infty
	\mbox{ as $\xi\to +\infty$} 
	\qquad \mbox{and} \qquad
	\int_{\xi_0}^\infty \frac{d\xi}{f(\xi)} < \infty
	\mbox{ for some $\xi_0>0$,}
\eas
then for some initial data the above solution must undergo some finite-time $L^\infty$ blow-up in the style 
described in ($\star$).\abs
\end{abstract}
\newpage
\section{Introduction}\label{intro}

Understanding the influence of nonlinearities has become a central theme in the literature on mathematical models 
for wave propagation during the past decades; especially for second-order problems, seminal works have revealed quite
subtle dependencies (see \cite{kenig_merle, merle_zaag, krieger_nakanishi_schlag_CMP2014, donninger, gazzola_squassina} and \cite{alinhac} for a small selection).
For physically relevant equations involving third-order time derivatives, however, corresponding knowledge to date 
appears to be much sparser. 
The present manuscript intends to undertake a first step toward the comprehension of blow-up phenomena
in a problem of this type.\Abs
{\bf The Jordan--Moore--Gibson--Thompson equation.} \quad
Motivated by the study of the propagation of nonlinear acoustic waves through thermally relaxing media, we consider the following
initial-boundary value problem:
\be{0} 
\lball
\tau u_{ttt} + \al u_{tt} = \beta \Del u_t + \gamma \Del u + \big(f(u)\big)_{tt},
\qquad & x\in\Om, \ t>0, \\[1mm]
u=0,
\qquad & x\in\pO, \ t>0, \\[1mm]
u(x,0)=u_0(x), \quad u_t(x,0)=u_1(x), \quad u_{tt}(x,0)=u_2(x),
\qquad & x\in\Om,
\ear
\ee
where the crucial nonlinearity suitably generalizes the prototypical setting determined by the choice
\be{exf}
f(u)=k u^2,  \qquad k >0.
\ee
The derivation of linear third-order equations,
\begin{equation} \label{MGT_eq}
	\tau u_{ttt} + \al u_{tt} = \beta \Del u_t + \gamma \Del u, 
\end{equation}
can be traced back to the work of Stokes~\cite{stokes1851xxxviii}. The model later appeared in the works of Moore and Gibson~\cite{moore1960propagation}, as well as Thompson~\cite{thompson}, and \eqref{MGT_eq} has since then
come to be known as the Moore--Gibson--Thompson (MGT) equation in the acoustic literature. 
It is derived by employing the Maxwell--Cattaneo relaxation of the Fourier heat flux law with the thermal relaxation time $\tau>0$ within the system of governing equations of sound motion as a means of avoiding the so-called \emph{paradox of infinite speed of propagation}.  
In the acoustic framework, $u$ represents the acoustic pressure, and if $\al=1$ is normalized for notational convenience, 
then the further system parameters are linked to the speed $c$ of sound in the medium 
and the so-called sound diffusivity $\delta$ via
the relations $\gamma=c^2$ and
\[
\beta = \delta+ \tau \gamma.
\]
The derivation of nonlinear acoustic models that incorporate thermal relaxation is due to Jordan~\cite{jordan2014second}, and they are thus commonly referred to as the Jordan--Moore--Gibson--Thompson (JMGT) equations. Two nonlinear versions of third-order models have become prominent in the literature, the second of which is covered by \eqref{0}. In~\cite[equation (72)]{jordan2014second}, the nonlinear third-order equation is given in the following form: 
\begin{equation} \label{JMGT_Kuznetsov}
	\tau \psi_{ttt} + \al \psi_{tt} = \beta \Del \psi_t + \gamma \Del u + \big(\tilde{\kappa} \psi_t^2+ \ell |\nabla \psi|^2 \big)_{t}.
\end{equation} 
Here, $\psi$ denotes the acoustic velocity potential, related to the acoustic pressure via $u = \varrho \psi_t$, where $\varrho$ is the medium density.  As \eqref{JMGT_Kuznetsov} can be understood as a hyperbolic version of the second-order strongly damped Kuznetsov equation~\cite{kuznetsov1971equations}: 
\begin{equation} \label{Kuznetsov}
	\al \psi_{tt} = \beta \Del \psi_t + \gamma \Del u + \big(\tilde{\kappa} \psi_t^2+ \ell |\nabla \psi|^2 \big)_{t},
\end{equation} 
formally obtained by setting $\tau$ to zero in \eqref{JMGT_Kuznetsov}, it is often referred to as the JMGT--Kuznetsov equation in mathematical works. When sound propagation is dominated by cumulative nonlinear effects, a simpler version of equation \eqref{JMGT_Kuznetsov} can be used. This will be the case when the propagation distance from the sound source is sufficiently large in terms of wavelengths, as discussed in~\cite{hamilton1998nonlinear}.
Approximating 
\begin{equation} \label{planewave_approx}
	\psi_t^2 \approx c^2 |\nabla \psi|^2;
\end{equation}
(see, e.g.~\cite{jordan2016survey, coulouvrat1992equations}), and then time-differentiating \eqref{JMGT_Kuznetsov} and expressing the resulting equation in terms of the acoustic pressure $u = \varrho \psi_t$ leads to the equation of present interest:  
\begin{equation} \label{JMGT_Westervelt}
	\tau u_{ttt} + \al u_{tt} = \beta \Del u_t + \gamma \Del u + \big(f(u)\big)_{tt}
\end{equation} 
with $f$ as in (\ref{exf}).
Since this popular form of the JMGT model can be seen as a hyperbolic version of the Westervelt equation~\cite{westervelt1963parametric} with strong damping,
\begin{equation} \label{Westervelt}
	\al u_{tt} = \beta \Del u_t + \gamma \Del u + \big(f(u)\big)_{tt},
\end{equation} 
it is also referred to as the JMGT--Westervelt equation.\abs
The study of nonlinear acoustic effects has become increasingly important with the rise of medical and industrial applications of ultrasound waves, which are inherently nonlinear because of their high frequencies and possible high intensities; see, e.g.,~the books~\cite{gallego2023power, szabo2004diagnostic} for an overview of these applications. While gradient blow-up is desirable in certain therapeutic ultrasound applications, such as shock-wave lithotripsy~\cite{cleveland2012physics}, it is usually avoided in diagnostic applications, such as imaging~\cite{szabo2004diagnostic}.  This observation highlights the need to investigate blow-up phenomena in the mathematical models of nonlinear ultrasonic propagation. \abs
However, at present, it seems unclear whether blow-up occurs in the third-order model \eqref{JMGT_Westervelt} and, if so, whether it is the gradient blow-up observed in the second-order model \eqref{Westervelt} when $\beta \rightarrow 0^+$. The primary purpose of the present work is thus to investigate and characterize blow-up phenomena in \eqref{0} with the nonlinearities such as \eqref{exf} relevant in nonlinear acoustic modeling, in a general setting of smooth bounded domains $\Omega \subset \R^n$ with $n \leq 3$.\Abs
{\bf State of art.} \quad
With regard to existence theories in ranges of solutions suitable far from singular behavior,
the mathematical research landscape on third-order acoustic equations is by now quite rich. The linear MGT equation has been extensively studied in the literature.  By the results of~\cite{kaltenbacher2011wellposedness}, the MGT model is known to generate a continuous semigroup which is exponentially stable when
\be{stable}
\alpha - \frac{\tau \gamma}{\beta}>0, \quad \gamma>0;
\ee
see~\cite[Theorem 1.2 and Theorem 1.4]{kaltenbacher2011wellposedness}.  The case $\alpha = \tau \gamma/\beta$ is referred to as the \emph{critical} or \emph{conservative} case where the generated semigroup is marginally stable.  The work of \cite{marchand2012abstract} complements \cite{kaltenbacher2011wellposedness} with theoretically precise exponential decay of the solutions using an operator-theoretic approach.  Regularity analysis in \cite{bucci2020regularity} has placed \eqref{MGT_eq} within a class of equations with memory to investigate the (interior) regularity of solutions. In a similar vein, a comparison to a viscoelastic model with memory can be found in~\cite{dell2017moore}. We mention that the influence of memory on the third-order model \eqref{MGT_eq} has also garnered a lot of interest in recent years; we refer to e.g., \cite{lasiecka2015moore, lasiecka2016moore, dell2016moore, alves2018moore} and the references contained therein.  \abs
The nonlinear problem \eqref{0}  with $f(u)= ku^2$ and $\alpha=1$ was first investigated in \cite{kaltenbacher2012well} on smooth bounded domains, where global well-posedness was established in the non-critical case $1-\frac{\tau \gamma}{\beta}>0$ assuming
\[
(u_0, u_1, u_2) \in  \left(\Htwoone \right) \times \Honezero \times \Ltwo
\]
sufficiently small;  see~\cite[Theorem 1.4]{kaltenbacher2012well}.  In~\cite{bongarti2020vanishing}, a rigorous relation to the Westervelt equation has been established through a singular limit analysis for the vanishing relaxation time $\tau$. A limiting local analysis for the vanishing sound diffusivity $\delta$ has been carried out in~\cite{kaltenbacher2021inviscid} under the assumption that  
\[
(u_0, u_1, u_2) \in  \left(\Htwoone \right) \times \left(\Htwoone \right) \times \Honezero.
\]
The well-posedness of a refined JMGT--Kuznetsov version of the model in \eqref{JMGT_Kuznetsov} has been investigated in \cite{kaltenbacher2019jordan}, uniformly in $\tau$, together with the accompanying singular limit analysis.   A global existence result and polynomial decay rates in time for a Cauchy problem in $\R^3$ for the same model have been obtained in~\cite{racke2021global}. A further improvement of these results in terms of data assumptions has been made in~\cite{said2022global}. \abs
A unidirectional approximation of \eqref{JMGT_Kuznetsov} by the hyperbolic Burgers' equation expressed in terms of $u_x$ has been used in~\cite[Section 4.4.1]{jordan2014second} to investigate blow-up phenomena. The results indicate finite time blow-up of $u_x$ for this approximate model for large data or a certain limiting value of the parameter $\tau$ specified in~\cite[equation (83)]{jordan2014second}. \abs
To the best of our knowledge, \change{the}
available results in the literature concerning blow-up directly for JMGT-type nonlinear equations 
concentrate on the mere detection of explosions, \change{and with the exception of~\cite{chen2022influence}}, in certain 
modified variants in which the forcing term $(f(u))_{tt}$ is replaced with different types of nonlinearities.
In \cite{chen2019blow}, a blow-up result for a Cauchy problem for the semilinear version of the equation with the right-hand side 
$|u_t|^p$  has been obtained  when the exponent satisfies $1< p \leq (n+1)/(n-1)$ for $n \geq 2$ and $p>1$ for $n=1$; see~\cite[Theorem 1.1]{chen2019blow}.  A version with the nonlinearity $|u|^p$ has been studied in~\cite{chen2019nonexistence} in the conservative and in~\cite{chen2021cauchy} in the dissipative case, and non-existence of global solutions has been proven under suitable assumptions on the exponent $p$. In~\cite{ming2021blow}, Cauchy problems for the versions of the JMGT equation with both such nonlinearities as well as their combination $|u|^p+|u_t|^p$ have been investigated further. \change{In~\cite[Section 2.4]{chen2022influence},  global nonexistence has been established for the Cauchy problem for the  JMGT-Kuznetsov and Westervelt equations in the conservative case.}
We underline, however, that even for variants of (\ref{0}) the knowledge on qualitative behavior of solutions
near blow-up seems fairly thin; specifically, the question how far solutions may cease to exist globally due to 
{\em derivative blow-up} has remained widely open even in simpler settings. 
\abs
We also note that,
in fact, the occurrence of
gradient blow-up has been numerically observed in the second-order Kuznetsov equation \eqref{Kuznetsov} and Westervelt equation \eqref{Westervelt} when $\beta \rightarrow 0^+$; see, e.g.,~\cite{walsh2007finite, kaltenbacher2007numerical, muhr2019}. In~\cite{dekkers2017cauchy}, 
the
Cauchy problem for the Kuznetsov equation in $\R^n$ has been studied, and a blow-up result established in the following sense:
\bas
\mbox{if $\tm<\infty$, \quad then \quad} \int_0^{\tm} \left\{ \|u_{tt}(\cdot,t)\|_{L^\infty(\R^n)}+\|\Delta u(\cdot,t)\|_{L^\infty(\R^n)}\right\}=\infty.
\eas
A further analytical study of the one-dimensional Kuznetsov equation, together with the investigation of the shock formation, can be found in~\cite{jordan2004analytical}. Having in mind the behavior of solutions to the second-order Kuznetsov model as well as the unidirectional Burgers' approximation used in \cite{jordan2004analytical}, a natural question is thus whether the same blow-up phenomena can be observed in \eqref{0} when $\tau>0$ and $\beta>0$.\Abs
{\bf Main results.} \quad
Our first objective is to make sure that unlike in the scenario from \cite{kaltenbacher2007numerical}, no type of
derivative blow-up can occur in (\ref{0}) throughout large classes of its ingredients, provided that with regard to their {first and}
third component, the initial data are assumed to be more regular than required in previous related literature on
local solvability (see Lemma \ref{lem10} below and~\cite{kaltenbacher2021inviscid}).
In fact, our main result in this direction applies to any physically meaningful
choice of the parameters $\tau,\beta,\gamma$ and $\al$, not necessarily linked through any relation of the form in (\ref{stable}),
and to widely arbitrary nonlinearities $f$:

\begin{theo}\label{theo13}
	Let $n\le 3$ and $\Om\subset\R^n$ be a bounded domain with smooth boundary, let $\tau>0,\beta>0,\gamma>0$ and $\al\in\R$,
	and suppose that
	\be{f}
	f\in  C^3(\R)
	\quad \mbox{is such that} \quad
	f(0)=0.
	\ee
	Then given any 
	\be{init}
	u_0\in {W^{4,2}(\Om)} \cap W_0^{1,2}(\Om),
	\ \
	u_1\in W^{2,2}(\Om) \cap W_0^{1,2}(\Om)
	\	\ \mbox{and} \
	u_2\in W^{2,2}(\Om) \cap W_0^{1,2}(\Om),
	\ee
	one can find $\tm\in (0,\infty]$ and a uniquely determined function $u$ such that
	\be{reg13}
	\lbal
	u\in C^0([0,\tm);W_0^{1,2}(\Om) \cap C^0(\bom)) \cap L^\infty_{loc}([0,\tm);W^{2,2}(\Om)), \\[1mm]
	u_t\in C^0([0,\tm);W_0^{1,2}(\Om) \cap C^0(\bom)) \cap L^\infty_{loc}([0,\tm);W^{2,2}(\Om)), \\[1mm]
	u_{tt}\in C^0([0,\tm);L^2(\Om)) \cap L^\infty_{loc}([0,\tm);W_0^{1,2}(\Om)), 
	\qquad \mbox{and} \\[1mm]
	u_{ttt} \in L^\infty_{loc}([0,\tm);L^2(\Om))
	\ear
	\ee
	that $u$ forms a strong solution of (\ref{0}) in the sense of Definition \ref{DS} below, and that
	\be{ext13}
	\mbox{if $\tm<\infty$, \quad then \quad}
	\limsup_{t\nearrow\tm} \|u(\cdot,t)\|_{L^\infty(\Om)} = \infty.
	\ee
\end{theo}
%
%
%
%
%
%
In a second step of our analysis, we will then proceed to identify conditions on the initial data under which
the above solutions must cease to exist in finite time, provided that $f$ is additionally assumed to be convex and
to satisfy mild assumptions on superlinear growth which particularly allow for the choice in (\ref{exf}).
In formulating our main result in this respect, as throughout the sequel we shall let, 
for a fixed bounded domain $\Om\subset\R^n$ with smooth boundary,
$(e_i)_{i\ge 1} \subset C^\infty(\bom)$ denote the system of Dirichlet eigenfunctions of $-\Del$, satisfying
\bas
\lball
-\Del e_i = \lam_i e_i,
\qquad & x\in\Om, \\[1mm]
e_i=0 & x\in\pO,
\ear
\eas
and normalized such that $\io e_i e_j = \del_{ij}$ for $i,j\in\N=\{1,2,3,...\}$, and that moreover
\bas
e_1>0
\qquad \mbox{in } \Om.
\eas
%
%
\begin{theo}\label{theo45}
	Let $n\le 3$ and $\Om\subset\R^n$ be a bounded domain with smooth boundary, and suppose that $f$ satisfies (\ref{f}) and
	is such that
	\be{f1}
	f''(\xi)\ge 0
	\qquad \mbox{for all } \xi\in\R,
	\ee
	that
	\be{f2}
	\frac{f(\xi)}{\xi} \to + \infty
	\qquad \mbox{as } \xi\to +\infty,
	\ee  
	and that
	\be{f3}
	\int_{\xi_0}^\infty \frac{d\xi}{f(\xi)} < \infty
	\qquad \mbox{for some } \xi_0>0.
	\ee
	Then for all $T_0>0$ there exists $K_0=K_0(T_0)>0$ such that if $u_0\in {W^{4,2}(\Om)}\cap W_0^{1,2}(\Om)$ satisfies
	\be{45.1}
	\io u_0 e_1 > K_0,
	\ee
	then one can find $K_1=K_1(u_0)>0$ in such a way that whenever $u_1\in W^{2,2}(\Om)\cap W_0^{1,2}(\Om)$ is such that
	\be{45.2}
	\io u_1 e_1 > K_1,
	\ee
	there is $K_2=K_2(u_1)>0$ such that if $u_2\in W^{2,2}(\Om)\cap W_0^{1,2}(\Om)$ has the property that
	\be{45.3}
	\io u_2 e_1 > K_2,
	\ee
	then the corresponding strong solution $u$ of (\ref{0}) blows up by time $T_0$ in the sense that in Theorem \ref{theo13} we
	have $\tm\le T_0$ and
	\be{45.4}
	\limsup_{t\nearrow \tm} \|u(\cdot,t)\|_{L^\infty(\Om)} = \infty.
	\ee
\end{theo}
~\abs
The rest of the exposition is organized as follows. In Section~\ref{Sec:local_existence}, local-in-time existence of strong solutions is shown together with an extensibility criterion, which forms the basis
of the exclusion of gradient blow-up phenomena.  Section~\ref{Sec:Excluding_grad_blowup} is then dedicated to the derivation of a refined extensibility criterion. The results of Section~\ref{sec:Uniqueness} complete the proof of Theorem~\ref{theo13} by establishing the uniqueness of a solution. Finally, Section~\ref{Sec:Blow-up} is devoted to the detection of finite-time $L^\infty$ blow-up and the proof of Theorem~\ref{theo45}.

\mysection{A preliminary statement on local existence and extensibility} \label{Sec:local_existence}
The purpose of this first step of our analysis will not only be 
to make sure by means of a fairly straightforward Galerkin-type approximation 
that for all suitably regular initial data, local-in-time strong solutions
always can be found; in its most substantial part, the main result of this section, as contained in Lemma \ref{lem8},
will also provide an extensibility criterion which, although yet being quite far from that in (\ref{ext13}), will form the basis
of our exclusion of gradient blow-up phenomena in Section \ref{Sec:Excluding_grad_blowup}.\abs
We already announce here that in line with precedent related works 
(cf.~, e.g., \cite{kaltenbacher2012well}),
developing this basic solution theory, and hence the analysis in this entire section, will only rely on the requirements that
\be{init0}
u_0\in W^{2,2}(\Om) \cap W_0^{1,2}(\Om),
\qquad
u_1\in W^{2,2}(\Om) \cap W_0^{1,2}(\Om)
\qquad \mbox{and} \qquad
u_2\in W_0^{1,2}(\Om),
\ee
which with respect to the assumptions on {$u_0$ and} $u_2$ are weaker than those in (\ref{init}).\abs
To begin with, let us specify the concept of solvability that will be pursued below.
\begin{defi}\label{DS}
	Let $n\ge 1$ and $\Om\subset\R^n$ be a bounded domain with smooth boundary, let $\tau>0,\beta>0,\gamma>0$, and $\al\in\R$,
	and let $f\in C^2(\R)$ and $T\in (0,\infty]$. 
	Then a function $u$ fulfilling
	\be{reg_S}
	\lbal
	u\in C^0([0,T);W_0^{1,2}(\Om) \cap C^0(\bom)) \cap L^\infty_{loc}([0,T);W^{2,2}(\Om)), \\[1mm]
	u_t\in C^0([0,T);W_0^{1,2}(\Om) \cap C^0(\bom)) \cap L^\infty_{loc}([0,T);W^{2,2}(\Om)), \\[1mm]
	u_{tt}\in C^0([0,T);L^2(\Om)) \cap L^\infty_{loc}([0,T);W_0^{1,2}(\Om)),
	\qquad \mbox{and} \\[1mm]
	u_{ttt} \in L^\infty_{loc}([0,T);L^2(\Om))
	\ear
	\ee
	will be called a {\em strong solution of (\ref{0})} in $\Om\times (0,T)$ if
	\bas
	u(\cdot,0)=u_0,
	\qquad
	u_t(\cdot,0)=u_1
	\qquad \mbox{and} \qquad
	u_{tt}(\cdot,0)=u_2,
	\eas
	and if
	
	\be{00}
	\tau u_{ttt} +{\al u_{tt}} = \beta \Del u_t + \gamma \Del u + f'(u) u_{tt} + f''(u) u_t^2
	\qquad \mbox{a.e.~in } \Om\times (0,T).
	\ee
\end{defi}
From now on, without explicitly mentioning this any further we shall suppose throughout that
$n\le 3$ and $\Om\subset\R^n$ is a smoothly bounded domain, that $\tau>0,\beta>0,\gamma>0$, and $\al\in\R$,
and that (\ref{f}) holds.\abs
Now an approximation of (\ref{0}) can be achieved by means of a straightforward Galerkin type discretization. To this end, we introduce
the approximations $u_{0N}, u_{1N}$ and $u_{2N}$ to given functions $u_0, u_1$ and $u_2$ fulfilling (\ref{init0}) by letting
\be{initN}
u_{kN}(x):=\sum_{i=1}^N \bigg\{ \io u_k e_i \bigg\} \cdot e_i(x),
\qquad x\in\bom,
\quad N\in\N,
\ee
for $k\in\{0,1,2\}$.
For convenience in notation in the next lemma and below, let us moreover introduce the finite-dimensional subspaces 
$V_N:={\rm span} \{ e_1,...,e_N \}$ of $W^{2,2}(\Om)\cap W_0^{1,2}(\Om)$ for $N\in\N$.

\begin{lem}\label{lem5}
	Suppose that (\ref{init0}) holds, 
	and for $k\in \{0,1,2\}$, let $(u_{kN})_{N\in\N}$ be as in (\ref{initN}).
	Then for each $N\in\N$ there exist $\tmN\in (0,\infty]$ as well as a function
	\be{5.1}
	u_N\in C^3([0,\tmN);V_N)
	\ee
	such that
	\be{5.2}
	u_N(\cdot,0)=u_{0N},
	\qquad 
	u_{Nt}(\cdot,0)=u_{1N}
	\qquad \mbox{and} \qquad
	u_{Ntt}(\cdot,0)=u_{2N},
	\ee
	so that for any $\vp\in C^0((0,\tmN);V_N)$ we have
	\be{5.3}
	\tau \io u_{Nttt}\vp + \al \io u_{Ntt}\vp
	= \beta \io \Del u_{Nt}\vp
	+ \gamma \io \Del u_N \vp
	+ \io \big(f(u_N)\big)_{tt} \vp
	\quad \mbox{for all } t\in (0,\tmN),
	\ee
	and that
	\be{5.4}
	\mbox{if $\tmN<\infty$, \quad then \quad}
	\limsup_{t\nearrow\tmN} \Big\{ \|u_N(\cdot,t)\|_{L^2(\Om)} + \|u_{Nt}(\cdot,t)\|_{L^2(\Om)} 
	+ \|u_{Ntt}(\cdot,t)\|_{L^2(\Om)} \Big\}
	= \infty.
	\ee
\end{lem}
\proof
In quite a standard manner, (\ref{5.2})--(\ref{5.3}) can be rewritten 
as an initial-value problem for a system of $N$ third-order ordinary differential
equations driven by locally Lipschitz continuous sources (cf.~\cite{temam2012infinite}). 
In view of the equivalence of norms on the finite-dimensional space $V_N$, the claim including the extensibility criterion
in (\ref{5.4}) thus immediately results from the Picard--Lindel\"of theorem.
\qedd 
The above approximation will turn out to cooperate well with one of the fundamental energy structures associated with the linear
Moore--Gibson--Thompson equation, at least when restricted to suitably small time intervals. 
In Lemma \ref{lem7} below, this will be seen to result from the following preliminary observation in this regard.
\begin{lem}\label{lem6}
	Assume (\ref{init0}), and let
	\be{B}
	B:={\frac{4\gamma^2}{\beta}}.
	\ee
	Then for each $N\in\N$, the function $\F_N \in C^1([0,\tmN))$ defined by setting
	\bea{FN}
	\F_N(t)
	&:=&\frac{\tau}{2} \io |\na u_{Ntt}(\cdot,t)|^2
	+ \frac{\beta}{2} \io |\Del u_{Nt}(\cdot,t)|^2 \nn\\
	& & + \gamma \io \Del u_N(\cdot,t) \Del u_{Nt}(\cdot,t)
	+ \frac{B}{2} \io |\Del u_N(\cdot,t)|^2,
	\qquad t\in [0,\tmN),
	\eea
	satisfies
	\be{6.1}
	\F_N(t) \ge \frac{\tau}{2} \io |\na u_{Ntt}(\cdot,t)|^2
	+ \frac{\beta}{4} \io |\Del u_{Nt}(\cdot,t)|^2
	+ \frac{B}{4} \io |\Del u_N(\cdot,t)|^2
	\qquad \mbox{for all } t\in (0,\tmN)
	\ee
	and
	\bea{6.2}
	\F_N'(t)
	&\le& \big\{ \|f'(u_N)\|_{L^\infty(\Om)} + |\al|+3\big\} \cdot \io |\na u_{Ntt}|^2
	+ (\gamma+1) \io |\Del u_{Nt}|^2
	+ \frac{B^2}{4} \io |\Del u_N|^2 \nn\\
	& & + \|f''(u_N)\|_{L^\infty(\Om)}^2 \|u_{Nt}\|_{L^\infty(\Om)}^2 \|\na u_{Nt}\|_{L^2(\Om)}^2 
	+ \frac{1}{4} \|f''(u_N)\|_{L^\infty(\Om)}^2 \|u_{Ntt}\|_{L^4(\Om)}^2 \|\na u_N\|_{L^4(\Om)}^2 \nn\\
	& & + \frac{1}{4} \|f'''(u_N)\|_{L^\infty(\Om)}^2 \|u_{Nt}\|_{L^\infty(\Om)}^4 \|\na u_N\|_{L^2(\Om)}^2
	\qquad \mbox{for all } t\in (0,\tmN).
	\eea
\end{lem}
\proof
Since
\bas
\bigg| \gamma \io \Del u_N \Del u_{Nt} \bigg| \le \frac{\beta}{4} \io |\Del u_{Nt}|^2
+ \frac{\gamma^2}{\beta} \io |\Del u_N|^2
\qquad \mbox{for all } t\in (0,\tmN)
\eas
by Young's inequality, the estimate in (\ref{6.1}) directly results from (\ref{FN}).\abs
To derive (\ref{6.2}), we note that according to our selection of $(e_i)_{i\in\N}$, for each $N\in\N$ the function
$\vp:=-\Del u_{Ntt}$ belongs to $C^0([0,\tmN);V_N)$, so that we may use (\ref{5.3}) and integrate by parts in a straightforward
manner to see that
for all $t\in (0,\tmN)$,
\bea{6.3}
& & \hs{-30mm}
\frac{\tau}{2} \frac{d}{dt} \io |\na u_{Ntt}|^2
+ \frac{\beta}{2} \frac{d}{dt} \io |\Del u_{Nt}|^2 \nn\\
&=& - \al \io |\na u_{Ntt}|^2
- \gamma { \io \Del u_N \Del u_{Ntt}}
+ \io \na \big( f(u_N) \big)_{tt} \cdot\na u_{Ntt} \nn\\
&=& - \al \io |\na u_{Ntt}|^2
- \gamma \frac{d}{dt} \io \Del u_N \Del u_{Nt}
+ \gamma \io |\Del u_{Nt}|^2 \nn\\
& & + \io f'(u_N) |\na u_{Ntt}|^2
+ \io f''(u_N) u_{Ntt} \na u_N\cdot\na u_{Ntt} \nn\\
& & + 2 \io f''(u_N) u_{Nt} \na u_{Nt}\cdot\na u_{Ntt} 
+ \io f'''(u_N) u_{Nt}^2 \na u_N \cdot\na u_{Ntt}.
\eea
Here, clearly,
\be{6.4}
- \al \io |\na u_{Ntt}|^2
+ \io f'(u_N) |\na u_{Ntt}|^2
\le \big\{ \|f'(u_N)\|_{L^\infty(\Om)} + |\al|\big\} \cdot \io |\na u_{Ntt}|^2
\qquad \mbox{for all } t\in (0,\tmN),
\ee
whereas Young's inequality and the Cauchy--Schwarz inequality ensure that
\bea{6.5}
&&\io f''(u_N) u_{Ntt} \na u_N\cdot\na u_{Ntt} \nn \\
&\le& \io |\na u_{Ntt}|^2
+ \frac{1}{4} \io f''^2(u_N) u_{Ntt}^2 |\na u_N|^2 \nn\\
&\le& \io |\na u_{Ntt}|^2
+ \frac{1}{4} \|f''(u_N)\|_{L^\infty(\Om)}^2 \|u_{Ntt}\|_{L^4(\Om)}^2 \|\na u_N\|_{L^4(\Om)}^2
\eea
and
\bea{6.6}
&&2 \io f''(u_N) u_{Nt} \na u_{Nt} \cdot \na u_{Ntt} \nn \\
&\le& \io |\na u_{Ntt}|^2
+ \io f''^2(u_N) u_{Nt}^2 |\na u_{Nt}|^2 \nn\\
&\le& \io |\na u_{Ntt}|^2
+ \|f''(u_N)\|_{L^\infty(\Om)}^2 \|u_{Nt}\|_{L^\infty(\Om)}^2 \|\na u_{Nt}\|_{L^2(\Om)}^2
\eea
as well as
\bea{6.7}
&&\io f'''(u_N) u_{Nt}^2 \na {u_N}\cdot\na u_{Ntt} \nn \\
&\le& \io |\na u_{Ntt}|^2
+ \frac{1}{4} \io f'''^2(u_N) u_{Nt}^4 |\na u_N|^2 \nn\\
&\le& \io |\na u_{Ntt}|^2
+ \frac{1}{4} \|f'''(u_N)\|_{L^\infty(\Om)}^2 \|u_{Nt}\|_{L^\infty(\Om)}^4 \|\na u_N\|_{L^2(\Om)}^2
\eea
for all $t\in (0,\tmN)$.
As furthermore, again by Young's inequality,
\be{6.8}
\frac{B}{2} \frac{d}{dt} \io |\Del u_N|^2
= B \io \Del {u_N} \Del u_{Nt}
\le \io |\Del u_{Nt}|^2
+ \frac{B^2}{4} \io |\Del u_N|^2
\qquad \mbox{for all } t\in (0,\tmN),
\ee
the inequality in (\ref{6.2}) results upon combining (\ref{6.3}) with (\ref{6.4})-(\ref{6.8}).
\qedd
In turning the inequality (\ref{6.2}) into expedient {\em a priori} estimates on time intervals of fixed $N$-independent length,
we shall now make use of our overall assumption $n\le 3$ for the first time.
\begin{lem}\label{lem7}
	For every $M>0$ there exist $T(M)>0$ and $C(M)>0$ such that whenever (\ref{init0}) holds with
	\be{7.1}
	\|\Del u_0\|_{L^2(\Om)}
	+ \|\Del u_1\|_{L^2(\Om)}
	+ \|\na u_2\|_{L^2(\Om)}
	\le M,
	\ee
	the solutions of (\ref{5.2})--(\ref{5.3}) have the properties that $\tmN\ge T(M)$ for all $N\in\N$, and that
	\be{7.2}
	\io |\na u_{Ntt}(\cdot,t)|^2 \le C(M)
	\qquad \mbox{for all $N\in\N$ and } t\in (0,T(M)),
	\ee  
	that
	\be{7.3}
	\io |\Del u_{Nt}(\cdot,t)|^2 \le C(M)
	\qquad \mbox{for all $N\in\N$ and } t\in (0,T(M)),
	\ee  
	that
	\be{7.4}
	\io |\Del u_N(\cdot,t)|^2 \le C(M)
	\qquad \mbox{for all $N\in\N$ and } t\in (0,T(M)),
	\ee  
	and {that}
	\be{7.5}
	\io u_{Nttt}^2(\cdot,t) \le C(M)
	\qquad \mbox{for all $N\in\N$ and } t\in (0,T(M)).
	\ee  
\end{lem}
\proof
As we are assuming that $n\le 3$, according to standard elliptic regularity theory and
the continuity of the embeddings $W^{2,2}(\Om) \hra L^\infty(\Om)$,
$W^{1,2}(\Om) \hra L^4(\Om)$ and $W^{2,2}(\Om) \hra W^{1,4}(\Om)$ we can find positive constants $c_1, c_2$ and $c_3$ such that
\be{7.6}
\|\vp\|_{L^\infty(\Om)} 
\le c_1\|\Del\vp\|_{L^2(\Om)}
\qquad \mbox{for all } \vp \in W^{2,2}(\Om) {\cap} W_0^{1,2}(\Om)
\ee
and
\be{7.7}
\|\vp\|_{L^4(\Om)}
\le c_2\|\na\vp\|_{L^2(\Om)}
\qquad \mbox{for all } \vp \in W_0^{1,2}(\Om)
\ee
as well as
\be{7.8}
\|\na\vp\|_{L^4(\Om)}
\le c_3\|\Del\vp\|_{L^2(\Om)}
\qquad \mbox{for all } \vp \in W^{2,2}(\Om) {\cap} W_0^{1,2}(\Om),
\ee
while a Poincar\'e inequality yields $c_4>0$ and $c_5>0$ fulfilling
\be{7.9}
\|\na\vp\|_{L^2(\Om)} \le {c_4}\|\Del\vp\|_{L^2(\Om)}
\qquad \mbox{for all } \vp \in W^{2,2}(\Om) {\cap} W_0^{1,2}(\Om)
\ee
and
\be{7.99}
\|\vp\|_{L^2(\Om)} \le c_5\|\na\vp\|_{L^2(\Om)}
\qquad \mbox{for all } \vp \in W_0^{1,2}(\Om).
\ee
Given $M>0$, with $B$ taken from (\ref{B}) we then let
\be{7.10}
c_6\equiv c_6(M):=\frac{\tau}{2} M^2 + \frac{\beta}{2} M^2 + \gamma M^2 + \frac{B}{2} M^2,
\ee
and abbreviating
\be{7.11}
c_7\equiv c_7(M):=c_1\cdot\Big( \frac{4(c_6+1)}{B}\Big)^\frac{1}{2}
\ee
we rely on the inclusion $f\in C^3(\R)$ to pick $c_8=c_8(M)>0$ such that
\be{7.12}
|f'(\xi)| + |f''(\xi)| + |f'''(\xi)| \le c_8
\qquad \mbox{for all } \xi\in [-c_7,c_7].
\ee
We thereupon introduce $c_i=c_i(M), i\in\{9,...,14\}$, by letting
\bas
c_9:=({R1, R2}{c_8}+|\al|+3) \cdot \frac{2(c_6+1)}{\tau},
\qquad
c_{10}:=(\gamma+1)\cdot \frac{4(c_6+1)}{\beta}
\qquad \mbox{and} \qquad
c_{11}:=\frac{B^2}{4} \cdot \frac{4(c_6+1)}{B}
\eas
as well as
\bas
c_{12}:=c_1^2 c_4^2 c_8^2 \cdot \Big(\frac{4(c_6+1)}{\beta}\Big)^2,
\qquad
c_{13}:=\frac{1}{4} c_2^2 c_3^2 c_8^2 \cdot \frac{2(c_6+1)}{\tau} \cdot \frac{4(c_6+1)}{B}
\eas
and
\bas
c_{14}:=\frac{1}{4} c_1^4 c_4^2 c_8^2 \cdot \Big(\frac{4(c_6+1)}{\beta}\Big)^2 \cdot \frac{4(c_6+1)}{B},
\eas
and writing
\be{7.188}
c_{15}\equiv c_{15}(M):=c_9+c_{10}+c_{11}+c_{12}+c_{13}+c_{14},
\ee
we define
\be{7.19}
T(M):=\frac{1}{2{c_{15}}}.
\ee
Henceforth assuming (\ref{init0}) and (\ref{7.1}), to show that then for each $N\in\N$ the quantity $\tmN$ from Lemma \ref{lem5}
and the function $\F_N$ in (\ref{FN}) satisfy $\tmN\ge T(M)$ and
\be{7.20}
\F_N(t) \le c_6+1
\qquad \mbox{for all } t\in (0,T(M)),
\ee
we set
\be{7.21}
S_N:=\Big\{ \wh{T}\in (0,\tmN) \ \Big| \ \F_N(t) \le c_6+1
\mbox{ for all } t\in (0,\wh{T}) \Big\},
\ee
and note that $S_N$ is not empty, and hence $T_N:=\sup S_N$ a well-defined element of $(0,\tmN] \subset (0,\infty]$, because
$\F_N$ is continuous and satisfies 
\bea{7.211}
\F_N(0)
&=& \frac{\tau}{2} \io |\na u_{2N}|^2
+ \frac{\beta}{2} \io |\Del u_{1N}|^2
+ \gamma \io \Del u_{0N} \Del u_{1N}
+ \frac{B}{2} \io |\Del u_{0N}|^2 \nn\\
&\le& \frac{\tau}{2} \io |\na u_{2N}|^2
+ \frac{\beta}{2} \io |\Del u_{1N}|^2
+ \gamma \|\Del u_{0N}\|_{L^2(\Om)} \|\Del u_{1N}\|_{L^2(\Om)}
+ \frac{B}{2} \io |\Del u_{0N}|^2 \nn\\
&\le& \frac{\tau}{2} \io |\na u_2|^2
+ \frac{\beta}{2} \io |\Del u_1|^2
+ \gamma \|\Del u_0\|_{L^2(\Om)} \|\Del u_1\|_{L^2(\Om)}
+ \frac{B}{2} \io |\Del u_0|^2 \nn\\
&\le& \frac{\tau}{2} \cdot M^2 + \frac{\beta}{2} \cdot M^2 + \gamma\cdot M \cdot M + \frac{B}{2} \cdot M^2
= c_6
\eea
due to (\ref{FN}), the Cauchy--Schwarz inequality and our assumption (\ref{7.1}) combined with evident nonexpansiveness properties
of the projection operation performed in (\ref{initN}).
Moreover, recalling (\ref{6.1}) we see that
\be{7.22}
\io |\na u_{Ntt}|^2 \le \frac{2(c_6+1)}{\tau},
\quad
\io |\Del u_{Nt}|^2 \le \frac{4(c_6+1)}{\beta}
\quad \mbox{and} \quad
\io |\Del u_N|^2 \le \frac{4(c_6+1)}{B}
\ee
for all $t \in (0, T_N)$, where the rightmost inequality together with (\ref{7.6}) and (\ref{7.11}) ensures that, in particular,
\be{7.23}
\|u_N\|_{L^\infty(\Om)} \le c_7
\qquad \mbox{for all } t\in (0,T_N).
\ee
Therefore, (\ref{7.12}) applies so as to show that on the right-hand side of (\ref{6.2}), again by (\ref{7.22}), we have
\bas
\big\{ \|f'(u_N)\|_{L^\infty(\Om)} + |\al|+3\big\} \cdot \io |\na u_{Ntt}|^2
\le (c_8+|\al|+3) \cdot \frac{2(c_6+1)}{\tau}
= c_9
\eas
and
\bas
(\gamma+1) \io |\Del u_{Nt}|^2
\le (\gamma+1) \cdot \frac{4(c_6+1)}{\beta} 
= c_{10}
\eas
as well as
\bas
\frac{B^2}{4} \io |\Del u_N|^2
\le \frac{B^2}{4} \cdot \frac{4(c_6+1)}{B}
= c_{11}
\eas
for all $t\in (0,T_N)$, while using (\ref{7.12}) in conjunction with (\ref{7.6})--(\ref{7.9}) we see that
\bas
\|f''(u_N)\|_{L^\infty(\Om)}^2 \|u_{Nt}\|_{L^\infty(\Om)}^2 \|\na u_{Nt}\|_{L^2(\Om)}^2
&\le& c_8^2 \cdot c_1^2 \|\Del u_{Nt}\|_{L^2(\Om)}^2 \cdot c_4^2 \|\Del u_{Nt}\|_{L^2(\Om)}^2 \\
&\le& c_1^2 c_4^2 c_8^2 \cdot \Big(\frac{4(c_6+1)}{\beta} \Big)^2
= c_{12}
\eas
and
\bas
\frac{1}{4} \|f''(u_N)\|_{L^\infty(\Om)}^2 \|u_{Ntt}\|_{L^4(\Om)}^2 \|\na u_N\|_{L^4(\Om)}^2
&\le& \frac{1}{4} c_8^2 \cdot c_2^2 \|\na u_{Ntt}\|_{L^2(\Om)}^2 \cdot c_3^2 \|\Del u_N\|_{L^2(\Om)}^2 \\
&\le& \frac{1}{4} c_2^2 c_3^2 c_8^2 \cdot \frac{2(c_6+1)}{\tau} \cdot \frac{4(c_6+1)}{B}
= c_{13}
\eas
and
\bas
\frac{1}{4} \|f'''(u_N)\|_{L^\infty(\Om)}^2 \|u_{Nt}\|_{L^\infty(\Om)}^4 \|\na u_N\|_{L^2(\Om)}^2
&\le& \frac{1}{4} c_8^2 \cdot c_1^4 \|\Del u_{Nt}\|_{L^2(\Om)}^4 \cdot c_4^2 \|\Del u_N\|_{L^2(\Om)}^2 \\
&\le& \frac{1}{4} c_1^4 c_4^2 c_8^2 \cdot \Big(\frac{4(c_6+1)}{\beta}\Big)^2 \cdot \frac{4(c_6+1)}{B}
= c_{14}
\eas
for all $t\in (0,T_N)$.
In view of (\ref{6.2}) and (\ref{7.188}), it thus follows that
\bas
\F_N'(t) \le c_{15}
\qquad \mbox{for all } t\in (0,T_N),
\eas
so that, once more by (\ref{7.211}),
\be{7.24}
\F_N(t) \le \F_N(0) + c_{15} t \le c_6 + c_{15} t
\qquad \mbox{for all } t\in (0,T_N).
\ee
Now if $T_N$ was smaller than $T(M)$, then, again by continuity of $\F_N$, either $\tmN\le T(M)$, or $\tmN>T(M)$ 
and $\F_N(T(M))=c_6+1$.
Here in light of (\ref{6.1}) the former alternative is clearly incompatible with (\ref{7.24}) 
due to the extensibility criterion in (\ref{5.4}); but also the latter possibility can readily be ruled out, for our definition
(\ref{7.19}) together with (\ref{7.24}) would show that in this case we should have
\bas
c_6+1 = \F_N(T(M)) \le c_6 + c_{15} T(M) \le c_6+\frac{1}{2},
\eas
which is absurd.\abs
The estimates in (\ref{7.2})--(\ref{7.4}) now directly result from (\ref{7.20}) when combined with (\ref{6.1}),
while (\ref{7.5}) can thereupon be concluded from (\ref{7.2})--(\ref{7.4}), because choosing $\vp:=u_{Nttt}$ in (\ref{5.3})
and employing the Cauchy--Schwarz inequality shows that
\bas
&&\tau \io u_{Nttt}^2 \\
&=& - \al \io u_{Ntt} u_{Nttt}
+ \beta \io \Del u_{Nt} u_{Nttt}
+ \gamma \io \Del u_N u_{Nttt} 
+ \io f'(u_N) u_{Ntt} u_{Nttt} \\& &\hspace*{9.2cm}
+ \io f''(u_N) u_{Nt}^2 u_{Nttt} \\
&\le& \Big\{ |\al|\cdot \|u_{Ntt}\|_{L^2(\Om)} + \beta \|\Del u_{Nt}\|_{L^2(\Om)} + \gamma \|\Del u_N\|_{L^2(\Om)} \\
& & \hs{10mm}
+ \|f'(u_N)\|_{L^\infty(\Om)} \|u_{Ntt}\|_{L^2(\Om)}
+ |\Om|^\frac{1}{2} \cdot \|f''(u_N)\|_{L^\infty(\Om)} \|u_{Nt}\|_{L^\infty(\Om)}^2 \Big\} \cdot \|u_{Nttt}\|_{L^2(\Om)}
\eas
and hence, by (\ref{7.99}) and (\ref{7.6}),
\bas
\tau \|u_{Nttt}\|_{L^2(\Om)}
&\le& |\al| c_5 \|\na u_{Ntt}\|_{L^2(\Om)}
+ \beta \|\Del u_{Nt}\|_{L^2(\Om)}
+ \gamma \|\Del u_N\|_{L^2(\Om)} \nn\\
& & + c_8 c_5 \|\na u_{Ntt}\|_{L^2(\Om)}
+ |\Om|^\frac{1}{2} c_8 c_1^2 \|\Del u_{Nt}\|_{L^2(\Om)}
\eas
for all $t\in (0,T(M))$ and each $N\in\N$.
\qedd
Through compactness properties thereby implied, the previous lemma leads us to the main result on this section:
\begin{lem}\label{lem8}
	If (\ref{init0}) holds, then there exist $\tm\in (0,\infty]$ and at least one strong solution $u$ of (\ref{0})
	in $\Om\times (0,\tm)$ which is such that 
	\be{ext8}
	\mbox{if $\tm<\infty$, \quad then \quad}
	\limsup_{t\nearrow\tm} \Big\{ \|\na u_{tt}(\cdot,t)\|_{L^2(\Om)} + \|\Del u_t(\cdot,t)\|_{L^2(\Om)}
	+ \|\Del u(\cdot,t)\|_{L^2(\Om)} \Big\} = \infty.
	\ee
\end{lem}
\proof
For fixed $M>0$ and any $(u_0,u_1,u_2)$ fulfilling (\ref{init0}) and (\ref{7.1}), Lemma \ref{lem7} {says} that with $T(M)$ and
$C(M)$ as provided there it follows that $\tmN\ge T(M)$ for all $N\in\N$, and that (\ref{7.2})-(\ref{7.5}) hold.
In particular, this means that
\bas
(u_N)_{N\in\N}
\mbox{ is bounded in } L^\infty((0,T(M));W^{2,2}(\Om) \cap W_0^{1,2}(\Om)),
\eas
that
\bas
(u_{Nt})_{N\in\N}
\mbox{ is bounded in } L^\infty((0,T(M));W^{2,2}(\Om) \cap W_0^{1,2}(\Om)),
\eas
that
\bas
(u_{Ntt})_{N\in\N}
\mbox{ is bounded in } L^\infty((0,T(M)); { W_0^{1,2}(\Om)}),
\eas
and that
\bas
(u_{Nttt})_{N\in\N}
\mbox{ is bounded in } L^\infty((0,T(M));L^2(\Om)),
\eas
whence due to an Aubin--Lions lemma (\cite{simon1986compact}) 
and the compactness of the embeddings $W^{2,2}(\Om) \hra W^{1,2}(\Om)$,
$W^{2,2}(\Om) \hra C^0(\bom)$ and $W^{1,2}(\Om) \hra L^2(\Om)$, we also know that
\bas
(u_N)_{N\in\N} 
\mbox{ is relatively compact in } C^0([0,T(M)];W_0^{1,2}(\Om) \cap C^0(\bom)),
\eas
that
\bas
(u_{Nt})_{N\in\N} 
\mbox{ is relatively compact in } C^0([0,T(M)];W_0^{1,2}(\Om) \cap C^0(\bom)),
\eas
and that
\bas
(u_{Ntt})_{N\in\N} 
\mbox{ is relatively compact in } C^0([0,T(M)];L^2(\Om)).
\eas
A straightforward extraction procedure thus yields $(N_j)_{j\in\N} \subset\N$ and a function 
\bas
u \in C^0([0,T(M)];W_0^{1,2}(\Om) \cap C^0(\bom)) \cap L^\infty((0,T(M));W^{2,2}(\Om))
\eas
satisfying
\bas
u_t \in C^0([0,T(M)];W_0^{1,2}(\Om) \cap C^0(\bom)) \cap L^\infty((0,T(M));W^{2,2}(\Om))
\eas
and
\bas
u_{tt} \in C^0([0,T(M)];L^2(\Om)) \cap L^\infty((0,T(M));W_0^{1,2}(\Om))
\eas
as well as
\bas
u_{ttt} \in L^\infty((0,T(M));L^2(\Om))
\eas
such that $N_j\to\infty$ as $j\to\infty$, and that
\begin{eqnarray}
	u_N \to u
	\quad & &
	\mbox{in } 
	C^0([0,T(M)];W_0^{1,2}(\Om) \cap C^0(\bom)), 
	\label{8.1} \\
	\Del u_N \wsto \Del u
	\quad & &
	\mbox{in } 
	L^\infty((0,T(M));L^2(\Om)),
	\label{8.2} \\
	u_{Nt} \to u_t
	\quad & &
	\mbox{in } 
	C^0([0,T(M)];W_0^{1,2}(\Om) \cap C^0(\bom)), 
	\label{8.3} \\
	\Del u_{Nt} \wsto \Del u_t
	\quad & &
	\mbox{in } 
	L^\infty((0,T(M));L^2(\Om)),
	\label{8.4} \\
	u_{Ntt} \to u_{tt}
	\quad & &
	\mbox{in } 
	C^0([0,T(M)];L^2(\Om)),
	\label{8.5} 
	\qquad \mbox{and} \\
	u_{Nttt} \wsto u_{ttt}
	\quad & & 
	\mbox{in }
	L^\infty((0,T(M));L^2(\Om))
	\label{8.6}
\end{eqnarray}
as $N=N_j\to\infty$.
To verify the claimed strong solution properties of $u$, we note that from (\ref{initN}) as well as (\ref{8.1}), (\ref{8.3})
and (\ref{8.5}) it directly follows that $u(\cdot,0)=\lim_{j\to\infty} u_{N_j}(\cdot,0)=\lim_{j\to\infty} u_{0 N_j}=u_0$
a.e.~in $\Om$, and that similarly also $u_t(\cdot,0)=u_1$ and $u_{tt}(\cdot,0)=u_2$ a.e.~in $\Om$.
Furthermore, for arbitrary $N_0\in\N$ and $\vp\in L^2((0,T(M));V_{N_0})$ it follows from (\ref{5.3}) that
\begin{equation}\label{8.7}
	\begin{aligned}
		&\tau \int_0^{T(M)} \io u_{Nttt} \vp
		+ \al \int_0^{T(M)} \io u_{Ntt} \vp \\
		=&\, \begin{multlined}[t] \beta \int_0^{T(M)} \io \Del u_{Nt} \vp + \gamma \int_0^{T(M)} \io \Del u_N \vp 
			+ \int_0^{T(M)} \io f'(u_N) u_{Ntt} \vp \\
			+ \int_0^{T(M)} \io f''(u_N) u_{Nt}^2 \vp \end{multlined}
	\end{aligned}
\end{equation}
for all $N\in\N$ such that $N\ge N_0$, so that since (\ref{8.1}) along with the continuity of $f'$ and $f''$
particularly ensures that $f'(u_N)\to f'({u})$ and 
$f''(u_N)\to f''(u)$ in $C^0(\bom\times [0,T(M)])$ as $N=N_j\to\infty$, and since (\ref{8.5}) and (\ref{8.3}) thus assert that
\bas
f'(u_N) u_{Ntt} \to f'(u) u_{tt}
\quad \mbox{and} \quad
f''(u_N) u_{Nt}^2 \to f''(u) u_t^2
\qquad \mbox{in } L^1(\Om\times (0,T(M))
\eas
as $N=N_j\to\infty$, from (\ref{8.7}) together with (\ref{8.6}), (\ref{8.5}), (\ref{8.4}) and (\ref{8.2}) we infer that
\bas
\tau \int_0^{T(M)} \io u_{ttt} \vp
+ \al \int_0^{T(M)} \io u_{tt} \vp
&=& \beta \int_0^{T(M)} \io \Del u_t \vp 
+ \gamma \int_0^{T(M)} \io \Del u \vp \\
& & + \int_0^{T(M)} \io f'(u) u_{tt} \vp
+ \int_0^{T(M)} \io f''(u) u_t^2 \vp.
\eas
Since $\bigcup_{N_0\in\N} V_{N_0}$ is dense in $L^2(\Om)$, this identity continues to hold actually for any 
$\vp\in L^2(\Om\times (0,T(M))$, so that indeed (\ref{00}) must hold.\abs
The blow-up criterion in (\ref{ext8}), finally, can be obtained by means of a standard extension argument, 
because our choice of $T(M)$ depended on $(u_0,u_1,u_2)$ only through the quantities appearing in (\ref{7.1}).
\qedd
\mysection{Excluding gradient blow-up} \label{Sec:Excluding_grad_blowup}
This section will now be devoted to the core part of our reasoning toward Theorem \ref{theo13}.
In particular, having Lemma \ref{lem8} at hand we may concentrate on the derivation of the refined extensibility
criterion \eqref{ext13} here, and this will be achieved by establishing bounds for the quantities
$\io |\na u_{tt}|^2$, $\io |\Del u_t|^2$, and $|\Del u|^2$ under the hypotheses that $\tm<\infty$, but that
$\sup_{t\in (0,\tm)} \|u(\cdot,t)\|_{L^\infty(\Om)}$ be finite.\abs
Of elementary nature but crucial importance in this regard will be the following observation related to certain integrated
versions of (\ref{0}).
{The use of these will later on, due to the appearance of $\Del u_0$ and 
	$u_2$ in (\ref{z1}), form the origin for our additional requirements in 
	(\ref{init}) beyond those from (\ref{init0}); indeed, Lemma~\ref{lem10} 
	will explicitly rely on $L^2$ bounds for the quantities $\Del z_1$ and 
	$\Del z_2$ determined by (\ref{z1}) and (\ref{z2}).}
\begin{lem}\label{lem9}
	Assume (\ref{init0}), and let $u$ be a strong solution of (\ref{0}) in $\Om\times (0,T)$ for some $T\in (0,\infty]$.
	Then writing
	\be{vw}
	v(x,t):=\int_0^t u(x,s) ds
	\quad \mbox{and} \quad
	w(x,t):=\int_0^t \int_0^s u(x,\sig) d\sig ds,
	\qquad
	(x,t)\in\Om\times (0,T),
	\ee
	as well as
	\be{z1}
	\begin{aligned}
		z_1(x):=&\, \tau u_2(x) + \al u_1(x) - \beta\Del u_0(x) - f'(u_0(x)) u_1(x)
		\quad x\in\Om,
	\end{aligned}
	\ee
	and
	\be{z2}
	\begin{aligned}
		z_2(x):=&\, \tau u_1(x) + \al u_0(x) - f(u_0(x))
		\quad x\in\Om,
	\end{aligned}
	\ee
	we have
	\be{01}
	\tau u_{tt} + \al u_t = \beta \Del u + \gamma \Del v + f'(u) u_t + z_1(x)
	\qquad \mbox{a.e.~in } \Om\times (0,T)
	\ee
	and
	\be{02}
	\tau u_t + \al u = \beta \Del v + \gamma \Del w + f(u) + tz_1(x) + z_2(x)
	\qquad \mbox{a.e.~in } \Om\times (0,T).
	\ee
\end{lem}
\proof 
Both identities immediately follow by integrating (\ref{00}) in time.
\qedd

For repeated later reference, let us separately state also the following collection of basic properties of resolvent-type
regularizations of functions from $L^2(\Om)$.
\begin{lem}\label{lem99}
	Let $A$ denote the realization of $-\Del$ in $L^2(\Om)$, with domain given by $D(A):=W^{2,2}(\Om) \cap W_0^{1,2}(\Om)$.
	Then for any $\eps>0$,
	\be{99.1}
	\Del (1+\eps A)^{-1} \vp = (1+\eps A)^{-1} \Del\vp
	\qquad \mbox{for all } \vp\in D(A)
	\ee
	and
	\be{99.2}
	\io \big| (1+\eps A)^{-1} \vp\big|^2
	\le \io \vp^2
	\qquad \mbox{for all } \vp\in L^2(\Om)
	\ee
	as well as
	\be{99.3}
	\io \big| \na \big\{ (1+\eps A)^{-1}\vp\big\}\big|^2
	\le \io |\na\vp|^2
	\qquad \mbox{for all } \vp\in W_0^{1,2}(\Om),
	\ee
	and for each fixed $\vp\in L^2(\Om)$, we moreover have
	\be{99.4}
	(1+\eps A)^{-1} \vp \to \vp
	\quad \mbox{in } L^2(\Om)
	\qquad \mbox{as } \eps\searrow 0.
	\ee
\end{lem}
{
	{
		\proof
		For $\vp\in L^2(\Om)$ and $\eps>0$, we let $\rho_\eps:=(1+\eps A)^{-1} \vp$, so that $\rho_\eps\in D(A)$ with
		\be{321}
		-\eps \Del\rho_\eps + \rho_\eps=\vp
		\qquad \mbox{a.e.~in } \Om.
		\ee
		If even $\vp\in D(A)$, then elliptic regularity theory (\cite[Theorem 1.17.2]{friedman2008partial}) asserts that 
		$\rho_\eps\in W^{4,2}(\Om)$, that $\Del\rho_\eps=\frac{1}{\eps} (\rho_\eps-\vp)$ belongs to $D(A)$, and that an application
		of $\Del=-A$ on both sides here yields 
		\bas
		-\eps \Del (\Del\rho_\eps) + (\Del\rho_\eps) = \Del\vp
		\qquad \mbox{a.e.~in } \Om,
		\eas
		meaning that, indeed, $\Del\rho_\eps\equiv \Del(1+\eps A)^{-1} \vp$ coincides with $(1+\eps A)^{-1} \Del\vp$.\abs
		To verify (\ref{99.2}), we only need to test (\ref{321}) against $\rho_\eps$ and employ Young's inequality, which namely
		shows that
		\bas
		\eps \io |\na\rho_\eps|^2 + \io \rho_\eps^2 = \io \vp\rho_\eps 
		\le \frac{1}{2} \io \rho_\eps^2 + \frac{1}{2} \io \vp^2,
		\eas
		and that thus $\frac{1}{2} \io \rho_\eps^2 \le \frac{1}{2} \io \vp^2$.\abs
		Similarly, if $\vp\in W_0^{1,2}(\Om)$, then multiplying (\ref{321}) by $-\Del\rho_\eps$ and integrating by parts we see that
		\bas
		\eps \io |\Del\rho_\eps|^2 + \io |\na\rho_\eps|^2 
		= \io \na\rho_\eps\cdot\na\vp
		\le \frac{1}{2} \io |\na\rho_\eps|^2
		+ \frac{1}{2} \io |\na\vp|^2,
		\eas
		from which we obtain that
		\be{322}
		\frac{1}{2} \io |\na\rho_\eps|^2 \le \frac{1}{2} \io |\na\vp|^2,
		\ee
		and thereby confirm (\ref{99.3}).\abs
		Finally, still for fixed $\vp\in W_0^{1,2}(\Om)$, (\ref{322}) and the Rellich embedding theorem particularly entail
		that the associated family $(\rho_\eps)_{\eps>0}$ is relatively compact in $L^2(\Om)$.
		To identify corresponding accumulation points along sequences of vanishing $\eps$, we only need to observe that,
		once more by (\ref{321}) and (\ref{322}), for any $\chi\in C_0^\infty(\Om)$ we have
		\bas
		\io \rho_\eps \chi = \io \vp\chi - \eps\io \na\rho_\eps\cdot\na\chi 
		\to \io \vp\chi
		\qquad \mbox{as } \eps\searrow 0.
		\eas
		We therefore conclude that for each $\vp\in W_0^{1,2}(\Om)$,
		\be{323}
		(1+\eps A)^{-1} \vp \to \vp
		\quad \mbox{in } L^2(\Om)
		\qquad \mbox{as } \eps\searrow 0.
		\ee
		To derive (\ref{99.4}) for arbitrary $\vp\in L^2(\Om)$, given any such $\vp$ and an arbitrary $\eta>0$ we choose
		$\vp_\eta\in W_0^{1,2}(\Om)$ such that $\|\vp-\vp_\eta\|_{L^2(\Om)} \le \frac{\eta}{3}$, and can then use (\ref{323})
		to find $\eps_\eta>0$ such that $\|(1+\eps A)^{-1} \vp_\eta - \vp_\eta\|_{L^2(\Om)} \le \frac{\eta}{3}$.
		Consequently, by linearity and (\ref{99.2}),
		\bas
		&&	\big\| (1+\eps A)^{-1} \vp - \vp \big\|_{L^2(\Om)}\\
		&\le& \big\| (1+\eps A)^{-1} (\vp_\eta-\vp)\big\|_{L^2(\Om)}
		+ \big\| (1+\eps A)^{-1} \vp_\eta - \vp_\eta \big\|_{L^2(\Om)}
		+ \|\vp_\eta-\vp\|_{L^2(\Om)} \\
		&\le& 2\|\vp_\eta-\vp\|_{L^2(\Om)}
		+ \big\| (1+\eps A)^{-1} \vp_\eta-\vp_\eta\big\|_{L^2(\Om)} \\
		&\le& \frac{2\eta}{3} + \frac{\eta}{3} = \eta
		\eas
		for all $\eps\in (0,\eps_\eta)$.
		\qed
	}
}
~\\[2mm]
The apparently most crucial step toward our derivation of Theorem \ref{theo13} can now be accomplished by means of a
testing procedure applied to the lifted equation in (\ref{02}), arranged in a suitably careful manner so as to circumvent
obstacles possibly going along with a lack of regularity.
We underline that from now on, the full regularity requirements in (\ref{init}) are made use of.
\begin{lem}\label{lem10}
	Beyond requiring (\ref{init0}), assume that
	\be{10.1}
	{u_0 \in W^{4,2}(\Om)}\quad \text{ and} \quad	u_2\in W^{2,2}(\Om),
	\ee
	and suppose that the strong solution $u$ of (\ref{0}) from Lemma \ref{lem8} has the properties that $\tm<\infty$, but that
	with some $K>0$ we have
	\be{10.2}
	\|u(\cdot,t)\|_{L^\infty(\Om)} \le K
	\qquad \mbox{for all } t\in (0,\tm).
	\ee
	Then there exists $C>0$ such that
	\be{10.3}
	\io |\Del u(\cdot,t)|^2 \le C
	\qquad \mbox{for a.e.~} t\in (0,\tm).
	\ee
\end{lem}
\proof
We let $A$ be as in Lemma \ref{lem99}, and for $\eps>0$ we let
\bea{10.4}
\F_\eps^{(1)}(t)
&:=& \frac{\tau}{2} \io \big| (1+\eps A)^{-1} \Del u(\cdot,t)\big|^2
+ \frac{\beta}{2} \io \big| \na \big\{ (1+\eps A)^{-1} \Del v(\cdot,t)\big\} \big|^2 \nn\\
& & + \gamma \io \na \big\{ (1+\eps A)^{-1} \Del v(\cdot,t)\big\} \cdot \na \big\{ (1+\eps A)^{-1} \Del w(\cdot,t)\big\}
\nn\\
& & + \frac{B}{2} \io \big|\na  \big\{ (1+\eps A)^{-1} \Del w(\cdot,t)\big\}\big|,
\qquad t\in [0,\tm),
\eea
where $v, w$ and $B$ are as in (\ref{vw}) and (\ref{B}).
Here, writing \[(1+\eps A)^{-1} \Del u = - \frac{1}{\eps} u + \frac{1}{\eps} (1+\eps A)^{-1} u,\] we see that since
$(1+\eps A)^{-1}$ acts as a bounded operator on $L^2(\Om)$, the inclusions $u\in C^0([0,\tm);L^2(\Om))$ and
$u_t\in C^0([0,\tm);L^2(\Om))$ particularly contained in (\ref{reg_S}) ensure that 
\[[0,\tm)\ni t \mapsto \io \big| (1+\eps A)^{-1} \Del u(\cdot,t)\big|^2\] is continuously differentiable with
\bas
\frac{\tau}{2} \frac{d}{dt} \io \big| (1+\eps A)^{-1} \Del u\big|^2
&=& \tau \io (1+\eps A)^{-1} \Del u \cdot (1+\eps A)^{-1} \Del u_t \\
&=& \tau \io (1+\eps A)^{-1} \Del u \cdot \Del \big\{ (1+\eps A)^{-1} u_t \big\},
\qquad \mbox{for all } t\in (0,\tm)
\eas
according to (\ref{99.1}).
Likewise, noting that 
\[\na \big\{ (1+\eps A)^{-1} \Del v \big\} \equiv -\frac{1}{\eps} \int_0^t \na u + \frac{1}{\eps} (1+\eps A)^{-1} \int_0^t \na u
\]
and
$\pa_t \na \big\{ (1+\eps A)^{-1} \Del v\big\} \equiv - \frac{1}{\eps} \na u + \frac{1}{\eps} (1+\eps A)^{-1} \na u$,
and similarly also $\na \big\{ (1+\eps A)^{-1} \Del w\big\}$, $\pa_t \big\{ (1+\eps A)^{-1} \Del w\big\}$, {and $\pa_t \nabla \big\{ (1+\eps A)^{-1} \Del w\big\}$}
all belong to $C^0([0,\tm);L^2(\Om;\R^n))$ by (\ref{reg_S}), we find that
\[[0,\tm) \ni t \mapsto \io \big| \na \big\{ (1+\eps A)^{-1} \Del v(\cdot,t)\big\} \big|^2,\]
\[[0,\tm) \ni t \mapsto \io \na \big\{ (1+\eps A)^{-1} \Del v(\cdot,t)\big\} 
\cdot \na \big\{ (1+\eps A)^{-1} \Del w(\cdot,t) \big\}\]
and
$[0,\tm) \ni t \mapsto \io \big| \na \big\{ (1+\eps A)^{-1} \Del w(\cdot,t)\big\} \big|^2$
lie in $C^1([0,\tm))$ and satisfy
\bas
\frac{\beta}{2} \frac{d}{dt} \io \big| \na \big\{ (1+\eps A)^{-1} \Del v\big\} \big|^2
&=& \beta \io \na \big\{ (1+\eps A)^{-1} \Del v\big\}\cdot \na \big\{ (1+\eps A)^{-1}\Del u \big\} \\
&=& - \beta \io (1+\eps A)^{-1} \Del u \cdot\Del \big\{ (1+\eps A)^{-1} \Del v \big\}
\qquad \mbox{for all } t\in (0,\tm)
\eas
and
\bas
& & \hs{-12mm}
\gamma\frac{d}{dt} \io \na \big\{ (1+\eps A)^{-1}\Del v\big\} \cdot\na \big\{ (1+\eps A)^{-1}\Del w\big\} \\
&=& \gamma \io \na \big\{ (1+\eps A)^{-1} \Del u\big\}\cdot \na \big\{ (1+\eps A)^{-1} \Del w\big\} 
+ \gamma \io \big| \na\big\{ (1+\eps A)^{-1}\Del v\big\} \big|^2 \\
&=& - \gamma \io (1+\eps A)^{-1} \Del u \cdot \Del \big\{ (1+\eps A)^{-1} {\Del w}\big\}
+ \gamma \io \big| \na\big\{ (1+\eps A)^{-1}\Del v\big\} \big|^2
\ \ \mbox{for all } t\in (0,\tm)
\eas
as well as
\bas
\frac{B}{2} \frac{d}{dt} \io \big| \na \big\{ (1+\eps A)^{-1} \Del w\big\}\big|^2
= B \io \na \big\{ (1+\eps A)^{-1} \Del v\big\} \cdot \na \big\{ (1+\eps A)^{-1} \Del w\big\}
\qquad \mbox{for all } t\in (0,\tm)
\eas
according to two integrations by parts.
In summary, we thus infer that $\F_\eps^{(1)}\in C^1([0,\tm))$, and that in line with (\ref{02}),
\bea{10.5}
\frac{d}{dt} \F_\eps^{(1)}(t)
&=& \io (1+\eps A)^{-1} \Del u \cdot \Del \Big\{ (1+\eps A)^{-1}  \big\{ \tau u_t - \beta\Del v - \gamma\Del w\big\} \nn\\
& & + \gamma \io \big|\na \big\{ (1+\eps A)^{-1} \Del v\big\}\big|^2 \nn\\
& & + B \io \na \big\{ (1+\eps A)^{-1} \Del v\big\}\cdot\na \big\{ (1+\eps A)^{-1} \Del w\big\} \nn\\
&=& \io (1+\eps A)^{-1} \Del u \cdot \Del \Big\{ (1+\eps A)^{-1} \big\{ - \al u + f(u) + tz_1 + z_2 \big\} \nn\\
& & + \gamma \io \big|\na \big\{ (1+\eps A)^{-1} \Del v\big\}\big|^2 \nn\\
& & + B \io \na \big\{ (1+\eps A)^{-1} \Del v\big\}\cdot\na \big\{ (1+\eps A)^{-1} \Del w\big\} \nn\\
&=& -\al \io \big| (1+\eps A)^{-1} \Del u\big|^2
+ \io (1+\eps A)^{-1}  \Del u \cdot \Del \big\{ (1+\eps A)^{-1}  f(u)\big\} \nn\\
& & + t\io (1+\eps A)^{-1} \Del u \cdot \Del \big\{ (1+\eps A)^{-1} z_1\big\}
+ \io (1+\eps A)^{-1} \Del u \cdot \Del \big\{ (1+\eps A)^{-1} z_2\big\} \nn\\
& & + \gamma \io \big|\na \big\{ (1+\eps A)^{-1} \Del v\big\}\big|^2 \nn\\
& & + B \io \na \big\{ (1+\eps A)^{-1} \Del v\big\}\cdot\na \big\{ (1+\eps A)^{-1} \Del w\big\}
\qquad \mbox{for all } t\in (0,\tm).
\eea
Now our definition of $B$ ensures that due to Young's inequality,
\be{10.6}
\F_\eps^{(1)}(t)
\ge \frac{\tau}{2} \io \big| (1+\eps A)^{-1} \Del u\big|^2
+ \frac{\beta}{4} \io \big|\na \big\{ (1+\eps A)^{-1} \Del v\big\}\big|^2
+ \frac{B}{4} \io \big|\na\big\{ (1+\eps A)^{-1} \Del w\big\} \big|^2
\ee
for all $t\in (0,\tm)$, 
so that, again by Young's inequality,
\bea{10.7}
& & \hs{-14mm}
- \al \io \big| (1+\eps A)^{-1} \Del u\big|^2
+ \gamma \io \big|\na \big\{ (1+\eps A)^{-1} \Del v\big\}\big|^2 
+ B \io \na \big\{ (1+\eps A)^{-1} \Del v\big\}\cdot\na \big\{ (1+\eps A)^{-1} \Del w\big\} \nn\\
&\le& |\al| \io \big| (1+\eps A)^{-1} \Del u\big|^2
+ (\gamma+1) \io \big|\na \big\{ (1+\eps A)^{-1} \Del v\big\} \big|^2
+ \frac{B^2}{4} \io \big| \na \big\{ (1+\eps A)^{-1} \Del w\big\}\big|^2 \nn\\
&\le& c_1 \F_\eps^{(1)}(t)
\qquad \mbox{for all } t\in (0,\tm),
\eea
where $c_1:=\frac{2|\al|}{\tau} + \frac{4(\gamma+1)}{\beta} + B$.\abs
Next, combining Young's inequality with (\ref{99.1}) and (\ref{99.2}) we see that
\bea{10.8}
& & \hs{-26mm}
t\io (1+\eps A)^{-1} \Del u \cdot \Del \big\{ (1+\eps A)^{-1} z_1\big\}
+ \io (1+\eps A)^{-1} \Del u \cdot \Del \big\{ (1+\eps A)^{-1} z_2\big\} \nn\\
&=& t \io (1+\eps A)^{-1} \Del u \cdot (1+\eps A)^{-1} \Del z_1
+ \io (1+\eps A)^{-1} \Del u \cdot (1+\eps A)^{-1}  \Del z_2 \\
&\le& \frac{\tau}{2} \io \big| (1+\eps A)^{-1} \Del u\big|^2 
+ \frac{t^2}{\tau} \io \big| (1+\eps A)^{-1} \Del z_1\big|^2
+ \frac{1}{\tau} \io \big| (1+\eps A)^{-1} \Del z_2\big|^2 \nn\\
&\le& \frac{\tau}{2} \io \big| (1+\eps A)^{-1} \Del u\big|^2
+ \frac{t^2}{\tau} {\io} |\Del z_1|^2
+ \frac{1}{\tau} \io |\Del z_2|^2 \nn\\
&\le& \F_\eps^{(1)}(t) + c_2
\qquad \mbox{for all } t\in (0,\tm),
\eea
where our assumption in (\ref{10.1}) together with { \eqref{7.20}}  and the hypothesis $\tm<\infty$ guarantees that 
{
	\begin{equation}
		\begin{aligned}
			c_2:=& \frac{\tm^2}{\tau} \io |\Del z_1|^2 + \frac{1}{\tau} \io |\Del z_2|^2 \\
			=&\, \begin{multlined}[t] \frac{\tm^2}{\tau} \io |\tau \Del  u_2 + \al \Del u_1 - \beta\Del^2 u_0 - \Del [f'(u_0) u_1] |^2 \\+ \frac{1}{\tau} \io |\tau \Del u_1 + \al \Del u_0 - \Del [ f(u_0)] |^2
			\end{multlined}
		\end{aligned}
	\end{equation}
	is finite.}\abs
To finally estimate the second summand on the right of (\ref{10.5}), we first proceed similarly to before to find that
\bea{10.9}
& & \hs{-30mm}
\io (1+\eps A)^{-1} \Del u\cdot\Del\big\{ (1+\eps A)^{-1}  f(u)\big\} \nn\\
&=& \io (1+\eps A)^{-1} \Del u \cdot (1+\eps A)^{-1} \Del f(u) \nn\\
&\le& \frac{\tau}{2} \io \big| (1+\eps A)^{-1} \Del u\big|^2
+ \frac{1}{2\tau} \io \big| (1+\eps A)^{-1} \Del f(u)\big|^2 \nn\\
&\le& \F_\eps^{(1)}(t)
+ \frac{1}{2\tau} \io |\Del f(u)|^2 \nn\\
&=& \F_\eps^{(1)}(t)
+ \frac{1}{2\tau} \io \big| f'(u) \Del u + f''(u) |\na u|^2\big|^2 \nn\\
&\le& \F_\eps^{(1)}(t)
+ \frac{1}{\tau} \io f'^2(u) |\Del u|^2
+ \frac{1}{\tau} \io f''^2(u) |\na u|^4
\qquad \mbox{for all } t\in (0,\tm),
\eea
and here we employ a Gagliardo--Nirenberg inequality, which together with elliptic regularity theory provides $c_3>0$ such that
\bas
\io |\na\vp|^4 
\le c_3\|\Del\vp\|_{L^2(\Om)}^2 \|\vp\|_{L^\infty(\Om)}^2
\qquad \mbox{for all } \vp\in W^{2,2}(\Om)\cap W_0^{1,2}(\Om).
\eas
Since our assumption in (\ref{10.2}) along with (\ref{f}) ensures the existence of $c_4>0$ and $c_5>0$ such that
$f'^2(u) \le c_4$ and $f''^2(u)\le c_5$ in $\Om\times (0,\tm)$, namely, this enables us to estimate
\bas
& & \hs{-30mm}
\frac{1}{\tau} \io f'^2(u) |\Del u|^2
+ \frac{1}{\tau} \io f''^2(u) |\na u|^4 \\
&\le& \frac{c_4}{\tau} \io |\Del u|^2
+ \frac{c_5}{\tau} \io |\na u|^4 \\
&\le& \frac{c_4}{\tau} \io |\Del u|^2
+ \frac{c_3 c_5}{\tau} \cdot \bigg\{ \io |\Del u|^2 \bigg\} \cdot \|u\|_{L^\infty(\Om)}^2 \\
&\le& c_6 \io |\Del u|^2
\qquad \mbox{for all } t\in (0,\tm)
\eas
with $c_6:=\frac{c_4}{\tau} + \frac{c_3 c_5 K^2}{\tau}$.\abs
In conjunction with (\ref{10.9}), (\ref{10.6}), (\ref{10.7}) and (\ref{10.5}), this shows that
\bas
\frac{d}{dt} \F_\eps^{(1)}(t)	
\le (c_1+2) \F_\eps^{(1)}(t)
+ c_6 \io |\Del u|^2
+ c_2
\qquad \mbox{for all } t\in (0,\tm),
\eas
which upon an integration in time yields the inequality
\bea{10.10}
\F_\eps^{(1)}(t)
&\le& \F_\eps^{(1)}(0) \cdot e^{(c_1+2)t}
+ c_6 \int_0^t e^{(c_1+2)(t-s)} \io |\Del u(\cdot,s)|^2 ds 
+ c_2 \int_0^t e^{({c_1+2})(t-s)} ds \nn\\
&\le& c_7 + c_8 \int_0^t \io |\Del u|^2
\qquad \mbox{for all } t\in (0,\tm)
\eea
with $c_7:=\big\{ \frac{\tau}{2} \io |\Del u_0|^2 + \frac{c_2}{c_1+2} \big\} \cdot e^{(c_1+2)\tm}$  and
$c_8:=c_6 e^{(c_1+2)\tm}$, because the identities $v(\cdot,0)\equiv 0$ and $w(\cdot,0)\equiv 0$ guarantee that
\bas
\F_\eps^{(1)}(0)
= \frac{\tau}{2} \io \big| (1+\eps A)^{-1} \Del u_0\big|^2
\le \frac{\tau}{2} \io |\Del u_0|^2
\eas
by (\ref{99.2}), and because
\bas
\int_0^t e^{(c_1+2)(t-s)} ds
= \frac{1}{c_1+2} (e^{(c_1+2)t}-1) \le \frac{1}{c_1+2} e^{(c_1+2)\tm}
\qquad \mbox{for all } t\in (0,\tm).
\eas
In view of (\ref{10.6}), from (\ref{10.10}) we particularly obtain that
\bas
\frac{\tau}{2} \io \big| (1+\eps A)^{-1} \Del u(\cdot,t)\big|^2 
\le c_7 + c_8 \int_0^t \io |\Del u|^2
\qquad \mbox{for all } t\in (0,\tm),
\eas
whence due to (\ref{99.4}) and the inclusion $\Del u\in L^\infty_{loc}([0,\tm);L^2(\Om))$ we infer upon taking $\eps\searrow 0$
that letting $y(t):=\io |\Del u(\cdot,t)|^2$, $t\in (0,\tm)$, defines a function $y\in L^\infty_{loc}([0,\tm))$ which
satisfies
\bas
\frac{\tau}{2} y(t) \le c_7 + c_8 \int_0^t y(s) ds
\qquad \mbox{for a.e.~} t\in (0,\tm).
\eas
An application of a Gronwall-type inequality (cf.~\cite[Section 1.3.6]{zheng2004nonlinear})
therefore reveals that
\bas
y(t) \le \frac{2c_7}{\tau} e^{\frac{2c_8}{\tau} t}
\qquad \mbox{for a.e.~} t\in (0,\tm),
\eas
which, again thanks to the finiteness of $\tm$, yields (\ref{10.3}).
\qedd
On the basis of the previous lemma, a second step in our bootstrap-like regularity analysis may operate on (\ref{01})
in a comparatively direct manner, without any need for additional regularization:
\begin{lem}\label{lem11}
	If (\ref{init}) holds, and if $\tm<\infty$ and (\ref{10.2}) is satisfied with some $K>0$, then
	there exists $C>0$ such that
	\be{11.1}
	\io |\na u_t(\cdot,t)|^2 \le C
	\qquad \mbox{for a.e.~} t\in (0,\tm).
	\ee
\end{lem}
\proof
Since (\ref{reg_S}) warrants that $\na u_t$ and $\na u_{tt}$ belong to $L^\infty_{loc}([0,\tm);L^2(\Om;\R^n))$, and that
$\Del u, \Del u_t, \Del v \equiv \int_0^t \Del u$ and $\Del v_t\equiv \Del u$ lie in $L^\infty_{loc}([0,\tm);L^2(\Om))$,
it follows that 
$[0,\tm)\ni t\mapsto \io |\na u_t(\cdot,t)|^2$,
$[0,\tm)\ni t\mapsto \io |\Del u(\cdot,t)|^2$,
$[0,\tm)\ni t\mapsto \io \Del u(\cdot,t) \Del v(\cdot,t)$
and
$[0,\tm)\ni t\mapsto \io |\Del v(\cdot,t)|^2$
are all elements of $W^{1,\infty}_{loc}([0,\tm))$,
whence again taking $B$ from (\ref{B}) we obtain that
\bas
\F^{(2)}(t)
&:=& \frac{\tau}{2} \io |\na u_t(\cdot,t)|^2
+ \frac{\beta}{2} \io |\Del u(\cdot,t)|^2 \\
& & + \gamma \io \Del u(\cdot,t) \Del v(\cdot,t)
+ \frac{B}{2} \io |\Del v(\cdot,t)|^2,
\qquad t\in [0,\tm),
\eas
defines a function $\F^{(2)} \in W^{1,\infty}_{loc}([0,\tm))$ fulfilling
\bea{11.11}
\frac{d}{dt} \F^{(2)}(t)
&=& \tau \io \na u_t\cdot\na u_{tt}
+ \beta \io \Del u \Del u_t \nn\\
& & + \gamma \io \Del u_t \Del v 
+ \gamma \io \Del u \Del v_t
+ B \io \Del v\Del v_t \nn\\
&=& \io \Del u_t \cdot \big\{ - \tau u_{tt} + \beta\Del u + \gamma \Del v \big\} \nn\\
& & + \gamma \io |\Del u|^2
+ B \io \Del u\Del v
\qquad \mbox{for a.e.~} t\in (0,\tm).
\eea
Since furthermore, by Young's inequality and (\ref{B}),
\be{11.2}
\F^{(2)}(t)
\ge \frac{\tau}{2} \io |\na u_t|^2
+ \frac{\beta}{4} \io |\Del u|^2
+ \frac{B}{4} \io |\Del v|^2
\qquad \mbox{for a.e.~} t\in (0,\tm),
\ee
we may here use Young's inequality to see that
\bea{11.3}
\gamma \io |\Del u|^2
+ B \io \Del u\Del v
&\le& (\gamma+1) \io |\Del u|^2
+ \frac{B^2}{4} \io |\Del v|^2 \nn\\
&\le& \Big\{ \frac{4(\gamma+1)}{\beta} + B\Big\} \cdot \F^{(2)}(t)
\qquad \mbox{for a.e.~} t\in (0,\tm),
\eea
while a combination of (\ref{01}) with two integrations by parts and Young's inequality shows that
\bea{11.4}
\io \Del u_t \cdot \big\{ -\tau u_{tt} + \beta\Del u + \gamma\Del v\big\}
&=& \io \Del u_t \cdot \big\{ \al u_t - f'(u) u_t - z_1\big\} \nn\\
&=& - \al \io |\na u_t|^2 \nn
+ \io \big\{ f'(u) \na u_t + f''(u) u_t\na u\big\} \cdot \na u_t \nn\\
& & + \io \na z_1\cdot\na u_t \nn\\
&\le& (|\al|+1) \io |\na u_t|^2 
+ \io f'(u) |\na u_t|^2 \nn\\
&&+ \io f''(u) u_t \na u\cdot\na u_t \nn\\
& & + \frac{1}{4} |\na z_1|^2
\ \qquad \mbox{for a.e.~} t\in (0,\tm).
\eea
Now according to our hypothesis that (\ref{10.2}) be valid, we can pick $c_1>0$ and $c_2>0$ such that $|f'(u)|\le c_1$ and
$|f''(u)|\le c_2$ in $\Om\times (0,\tm)$, and therefore, by another application of Young's inequality, we see that
\bea{11.5}
& & \hs{-30mm}
\io f'(u) |\na u_t|^2
+ \io f''(u) u_t \na u\cdot\na u_t \nn\\
&\le& c_1 \io |\na u_t|^2
+ c_2 \io |u_t| \cdot |\na u| \cdot |\na u_t| \nn\\
&\le& (c_1+1) \io |\na u_t|^2
+ \frac{c_2^2}{4} \io u_t^2 |\na u|^2
\qquad \mbox{for a.e.~} t\in (0,\tm).
\eea
But if in line with elliptic regularity and the continuity of the embeddings $W^{1,2}(\Om)\hra L^4(\Om)$ and
$W^{2,2}(\Om)\hra W^{1,4}(\Om)$, as asserted by the inequality $n\le 3$, we fix $c_3>0$ and $c_4>0$ such that
\bas
\|\vp\|_{L^4(\Om)} \le c_3\|\na\vp\|_{L^2(\Om)}
\qquad \mbox{for all } \vp\in W_0^{1,2}(\Om)
\eas
and
\bas
\|\na\vp\|_{L^4(\Om)} \le c_4\|\Del\vp\|_{L^2(\Om)}
\qquad \mbox{for all } \vp\in W^{2,2}(\Om)\cap W_0^{1,2}(\Om),
\eas
we may here invoke the Cauchy--Schwarz inequality to estimate
\bas
\io u_t^2 |\na u|^2
\le { \|u_t\|_{L^4(\Om)}^2} \|\na u\|_{L^4(\Om)}^2
\le c_3^2 c_4^2 \|\na u_t\|_{L^2(\Om)}^2 \|\Del u\|_{L^2(\Om)}^2
\qquad \mbox{for a.e.~} t\in (0,\tm),
\eas
whence an application of Lemma \ref{lem10} yields $c_5>0$ fulfilling
\bas
\frac{c_2^2}{4} \io u_t^2 |\na u|^2
\le c_5 \io |\na u_t|^2
\qquad \mbox{for a.e.~} t\in (0,\tm).
\eas
From (\ref{11.4}), (\ref{11.5}) and (\ref{11.2}) we thus obtain that
\bas
& & \hs{-30mm}
\io \Del u_t\cdot\big\{ -\tau u_{tt} + \beta\Del u +\gamma\Del v\big\} \\
&\le& (|\al|+c_1+c_5+2) \io |\na u_t|^2
+ \frac{1}{4} \io |\na z_1|^2 \\
&\le& \frac{2(|\al|+c_1+c_5+2)}{\tau} \cdot \F^{(2)}(t)
+ \frac{1}{4} \io |\na z_1|^2
\qquad \mbox{for a.e.~} t\in (0,\tm),
\eas
and that therefore, by (\ref{11.11}) and (\ref{11.3}),
\bas
\frac{d}{dt} \F^{(2)}(t)
\le \Big\{ \frac{4(\gamma+1)}{\beta} + B + \frac{2(|\al|+c_1+c_5+2)}{\tau}\Big\} \cdot \F^{(2)}(t)
+ \frac{1}{4} \io |\na z_1|^2
\qquad \mbox{for a.e.~} t\in (0,\tm).
\eas
As $z_1 = \tau u_2+ \al u_1 - \beta\Del u_0 - f'(u_0) u_1\in W^{1,2}(\Om)$ by (\ref{init}), a Gronwall lemma guarantees boundedness of $\F^{(2)}$ on $(0,\tm)$, so that the
claim results upon once more utilizing (\ref{11.2}).
\qedd
The decisive third-order quantities in (\ref{ext8}), finally, can be controlled by using the original equation (\ref{00}) 
in the course of an analysis of an energy functional resembling that in (\ref{FN}).
Again, potentially insufficient regularity features of strong solutions suggest to perform some suitable regularization
in the argument.
\begin{lem}\label{lem12}
	Assume (\ref{init}), and suppose that $\tm<\infty$, and that (\ref{10.2}) holds with some $K>0$.
	Then there exists $C>0$ such that
	\be{12.1}
	\io |\na u_{tt}(\cdot,t)|^2 \le C
	\qquad \mbox{for a.e.~} t\in (0,\tm).
	\ee
	and
	\be{12.2}
	\io |\Del u_t(\cdot,t)|^2 \le C
	\qquad \mbox{for a.e.~} t\in (0,\tm).
	\ee
\end{lem}
\proof
For $\eps>0$, thanks to the boundedness of the operators $\na (1+\eps A)^{-1}$ and $\Del (1+\eps A)^{-1}$ on $L^2(\Om)$
it follows from (\ref{reg_S}) that
\[
\Big\{ \na \big\{ (1+\eps A)^{-1} u_{tt}\big\} \, , \, \na \big\{ (1+\eps A)^{-1} u_{ttt}\big\}\Big\} \subset
L^\infty_{loc}[0,\tm);L^2(\Om;\R^n)),
\]
and that
$\Big\{ \Del \big\{ (1+\eps A)^{-1} u\big\} \, , \, \Del \big\{ (1+\eps A)^{-1} u_t \big\} \, , \, 
\Del \big\{ (1+\eps A)^{-1} u_{tt}\big\}\Big\} \subset
L^\infty_{loc}([0,\tm);L^2(\Om))$.
Therefore, letting $B$ be as in (\ref{B}) and
\bea{12.22}
\F_\eps^{(3)}(t)
&:=& \frac{\tau}{2} \io \big| \na\big\{ (1+\eps A)^{-1} u_{tt}(\cdot,t)\big\} \big|^2
+ \frac{\beta}{2} \io \big| \Del \big\{ (1+\eps A)^{-1} u_t(\cdot,t)\big\} \big|^2 \nn\\
& & + \gamma \io \Del \big\{ (1+\eps A)^{-1}  u(\cdot,t)\big\} \cdot \Del \big\{ (1+\eps A)^{-1} u_t(\cdot,t)\big\} \nn\\
& & + \frac{B}{2} \io \big| \Del \big\{ (1+\eps A)^{-1} u(\cdot,t)\big\}\big|^2,
\qquad t\in [0,\tm),
\eea
we obtain a function $\F_\eps^{(3)} \in W^{1,\infty}_{loc}([0,\tm))$ which, besides satisfying
\begin{equation}\label{12.3}
	\begin{aligned}
		\F_\eps^{(3)}(t) 
		\ge\, \begin{multlined}[t]\frac{\tau}{2} \io \big|\na \big\{ (1+\eps A)^{-1}  u_{tt}\big\}\big|^2
			+ \frac{\beta}{4} \io \big| \Del \big\{ (1+\eps A)^{-1} u_t\big\}\big|^2
			\\+ \frac{B}{4} \io \big| \Del \big\{ (1+\eps A)^{-1} u\big\}\big|^2 \end{multlined}
	\end{aligned}
\end{equation}
for a.e.~$t\in (0,\tm)$ by Young's inequality and (\ref{B}), has the property that
\bea{12.4}
\frac{d}{dt} \F_\eps^{(3)}(t)
&=& \tau \io \na \big\{ (1+\eps A)^{-1} u_{tt}\big\} \cdot\na\big\{ (1+\eps A)^{-1}  u_{ttt}\big\} \nn\\
& & + \beta \io \Del \big\{ (1+\eps A)^{-1}  u_t\big\} \cdot \Del \big\{ (1+\eps A)^{-1} u_{tt}\big\} \nn\\
& & + \gamma \io \Del \big\{ (1+\eps A)^{-1} u\big\} \cdot \Del \big\{ (1+\eps A)^{-1}  u_{tt}\big\}
+ \gamma \io \big|\Del \big\{ (1+\eps A)^{-1}  u_t\big\} \big|^2 \nn\\
& & + B \io {\Delta} \big\{ (1+\eps A)^{-1} u \big\} \cdot \Del \big\{ (1+\eps A)^{-1}  u_t\big\}
\eea
for a.e.~$t\in (0,\tm)$.
Here, twice integrating by parts and using (\ref{00}) shows that {due} to (\ref{99.1}) and Young's inequality,
\bea{12.5}
& & \hs{-22mm}
\tau \io \na \big\{ (1+\eps A)^{-1} u_{tt}\big\} \cdot\na\big\{ (1+\eps A)^{-1}  u_{ttt}\big\}
+ \beta \io \Del \big\{ (1+\eps A)^{-1}  u_t\big\} \cdot \Del \big\{ (1+\eps A)^{-1} u_{tt}\big\} \nn\\
& & \hs{-22mm}
+ \gamma \io \Del \big\{ (1+\eps A)^{-1} u\big\} \cdot \Del \big\{ (1+\eps A)^{-1}  u_{tt}\big\} \nn\\
&=& \io \na\big\{ (1+\eps A)^{-1} u_{tt}\big\} \cdot\na \Big\{ (1+\eps A)^{-1} 
\big\{ \tau u_{ttt} - \beta \Del u_t - \gamma \Del u \big\} \Big\} \nn\\
&=& - \al \io \big| \na \big\{ (1+\eps A)^{-1}  u_{tt}\big\}\big|^2 \nn\\
& & + \io \na \big\{ (1+\eps A)^{-1} u_{tt}\big\}\cdot\na \Big\{ (1+\eps A)^{-1}  \big\{ f'(u) u_{tt}\big\} \Big\} \nn\\
& & + \io \na \big\{ (1+\eps A)^{-1} u_{tt}\big\} \cdot\na \Big\{ (1+\eps A)^{-1} \big\{ f''(u) u_t^2 \big\} \Big\} \nn\\
&\le& (|\al|+1) \io \big|\na \big\{ (1+\eps A)^{-1} u_{tt}\big\}\big|^2 \nn\\
& & + \frac{1}{2} \io \Big|\na \Big\{ (1+\eps A)^{-1}  \big\{ f'(u) u_{tt} \big\} \Big\} \Big|^2 \nn\\
& & + \frac{1}{2} \io \Big|\na\Big\{ (1+\eps A)^{-1} \big\{ f''(u) u_t^2 \big\} \Big\} \Big|^2
\qquad \mbox{for a.e.~} t\in (0,\tm).
\eea
To estimate the two latter integrals, we combine (\ref{10.2}) with (\ref{f}) in choosing $c_1>0, c_2>0$ and $c_3>0$ such that
$|f'(u)|\le c_1, |f''(u)|\le c_2$ and $|f'''(u)|\le c_3$ on $\Om\times (0,\tm)$, and employ (\ref{99.3}) in verifying that since
$f'(u(\cdot,t))u_{tt}(\cdot,t)$ and $f''(u(\cdot,t)) u_t^2(\cdot,t)$ clearly belong to $W_0^{1,2}(\Om)$ for a.e~$t\in (0,\tm)$
according to (\ref{reg_S}), we have
\bea{12.6}
\frac{1}{2} \io \Big|\na \Big\{ (1+\eps A)^{-1}  \big\{ f'(u) u_{tt} \big\} \Big\} \Big|^2 	
&\le& \frac{1}{2} \io \big| \na \big\{ f'(u) u_{tt}\big\} \big|^2 \nn\\
&=& \frac{1}{2} \io \big| f'(u) \na u_{tt} + f''(u) u_{tt} \na u \big|^2 \nn\\
&\le& \io f'^2(u) |\na u_{tt}|^2
+ \io f''^2(u) u_{tt}^2 |\na u|^2 \nn\\
&\le& c_1^2 \io |\na u_{tt}|^2
+ c_2^2 \io u_{tt}^2 |\na u|^2
\eea
and
\bea{12.7}
\frac{1}{2} \io \Big|\na\Big\{ (1+\eps A)^{-1} \big\{ f''(u) u_t^2 \big\} \Big\} \Big|^2 	
&\le& \frac{1}{2} \io \big| \na \big\{ f''(u) u_t^2\big\} \big|^2 \nn\\
&=& \frac{1}{2} \io \big| 2f''(u) u_t \na u_t + f'''(u) u_t^2 \na u\big|^2 \nn\\
&\le& 4 \io f''^2(u) u_t^2 |\na u_t|^2
+ \io f'''^2(u) u_t^4 |\na u|^2 \nn\\
&\le& 4c_2^2 \io u_t^2 |\na u_t|^2
+ c_3^2 \io u_t^4 |\na u|^2
\eea
for a.e.~$t\in (0,\tm)$.\abs
Now according to elliptic regularity theory and the continuity of all the embeddings $W^{1,2}(\Om)\hra L^3(\Om)$,
$W^{2,2}(\Om)\hra W^{1,6}(\Om), W^{2,2}(\Om) \hra L^\infty(\Om)$ and $W^{1,2}(\Om)\hra L^6(\Om)$,
we can pick positive constants $c_4, c_5, c_6$ and $c_7$ such that
\be{12.8}
\|\vp\|_{L^3(\Om)} \le c_4 \|\na\vp\|_{L^2(\Om)}
\qquad \mbox{for all } \vp\in W_0^{1,2}(\Om),
\ee
that
\be{12.9}
\|\na\vp\|_{L^6(\Om)} \le c_5\|\Del\vp\|_{L^2(\Om)}
\qquad \mbox{for all } \vp\in W^{2,2}(\Om) \cap W_0^{1,2}(\Om),
\ee
and that
\be{12.10}
\|\vp\|_{L^\infty(\Om)} \le c_6 \|\Del\vp\|_{L^2(\Om)}
\qquad \mbox{for all } \vp\in W^{2,2}(\Om) \cap W_0^{1,2}(\Om)
\ee
as well as
\be{12.11}	
\|\vp\|_{L^6(\Om)} \le c_7\|\na\vp\|_{L^2(\Om)}
\qquad \mbox{for all } \vp\in W_0^{1,2}(\Om),
\ee
whereas Lemma \ref{lem10} and Lemma \ref{lem11} show that due to (\ref{10.2}) there exist $c_8>0$ and $c_9>0$ such that
\be{12.12}
\|\Del u\|_{L^2(\Om)} \le c_8
\qquad \mbox{for a.e.~} t\in (0,\tm)
\ee
and
\be{12.13}
\|\na u_t\|_{L^2(\Om)} \le c_9
\qquad \mbox{for a.e.~} t\in (0,\tm).
\ee
On the right-hand side of (\ref{12.6}), on the basis of the H\"older inequality, (\ref{12.8}), (\ref{12.9}) and (\ref{12.12})
we thus obtain that
\bas
c_2^2 \io u_{tt}^2 |\na u|^2
&\le& c_2^2 \|u_{tt}\|_{L^3(\Om)}^2 \|\na u\|_{L^6(\Om)}^2 \\
&\le& c_2^2 c_4^2 c_5^2 \|\na u_{tt}\|_{L^2(\Om)}^2{ \|\Del u\|_{L^2(\Om)}^2} \\
&\le& c_2^2 c_4^2 c_5^2 c_8^2 \io |\na u_{tt}|^2
\qquad \mbox{for a.e.~} t\in (0,\tm),
\eas
while (\ref{12.10}) together with (\ref{12.13}) shows that
\bas
4c_2^2 \io u_t^2 |\na u_t|^2
&\le& 4c_2^2 \|u_t\|_{L^\infty(\Om)}^2 \|\na u_t\|_{L^2(\Om)}^2 \\
&\le& 4c_2^2 c_6^2 \|\Del u_t\|_{L^2(\Om)}^2 \|\na u_t\|_{L^2(\Om)}^2 \\
&\le& 4c_2^2 c_6^2 c_9^2 \io |\Del u_t|^2
\qquad \mbox{for a.e.~} t\in (0,\tm).
\eas
Furthermore, using the H\"older inequality along with (\ref{12.11}), (\ref{12.9}), (\ref{12.13}) and (\ref{12.12}) we find that
\bas
c_3^2 \io u_t^4 |\na u|^2
&\le& c_3^2 { \|u_t\|_{L^6(\Om)}^4} \|\na u\|_{L^6(\Om)}^2 \\
&\le& c_3^2 c_7^2 c_5^2 \|\na u_t\|_{L^2(\Om)}^4 \|\Del u\|_{L^2(\Om)}^2 \\
&\le& c_3^2 c_7^2 c_5^2 c_9^4 c_8^2
\qquad \mbox{for a.e.~} t\in (0,\tm),
\eas
whence in view of (\ref{10.3}) we infer from (\ref{12.5}) that
\bea{12.14}
& & \hs{-22mm}
\tau \io \na \big\{ (1+\eps A)^{-1} u_{tt}\big\} \cdot\na\big\{ (1+\eps A)^{-1}  u_{ttt}\big\}
+ \beta \io \Del \big\{ (1+\eps A)^{-1}  u_t\big\} \cdot \Del \big\{ (1+\eps A)^{-1} u_{tt}\big\} \nn\\
& & \hs{-22mm}
+ \gamma \io \Del \big\{ (1+\eps A)^{-1} u\big\} \cdot \Del \big\{ (1+\eps A)^{-1}  u_{tt}\big\} \nn\\
&\le& c_{10} \io |\na u_{tt}|^2 + c_{11} \io |\Del u_t|^2 + c_{12}
\qquad \mbox{for a.e.~} t\in (0,\tm)
\eea
with $c_{10}:=|\al|+1+c_1^2+c_2^2 c_4^2 c_5^2 c_8^2$, $c_{11}:=4c_2^2 c_6^2 c_9^2$ and
$c_{12}:=c_3^2 c_5^2 c_7^2 c_8^2 c_9^4$.
As the two last summands in (\ref{12.4}) can be controlled by means of Young's inequality and (\ref{12.3}) according to
\bas
& & \hs{-20mm}
\gamma \io \big| \Del \big\{ (1+\eps A)^{-1} u_t\big\} \big|^2
+ B \io \Del \big\{ (1+\eps A)^{-1} u \big\} \cdot \Del \big\{ (1+\eps A)^{-1} u_t \big\} \\
&\le& (\gamma+1) \io \big| \Del \big\{ (1+\eps A)^{-1} u_t\big\} \big|^2
+ \frac{B^2}{4} \io \big| \Del \big\{ (1+\eps A)^{-1} u\big\} \big|^2 \\
&\le& c_{13} \F_\eps^{(3)}(t)
\qquad \mbox{for a.e.~} t\in (0,\tm)
\eas
with $c_{13} :=\frac{4(\gamma+1)}{\beta}+B$, from (\ref{12.14}) and (\ref{12.3}) we all in all conclude that
\bas
\frac{d}{dt} \F_\eps^{(3)}(t) \le c_{13} \F_\eps^{(3)}(t)
+ c_{10} \io |\na u_{tt}|^2
+ c_{11} \io |\Del u_t|^2 + c_{12}
\qquad \mbox{for a.e.~} t\in (0,\tm)
\eas
and hence
\bea{12.15}
\F_\eps^{(3)}(t)
&\le& \F_\eps^{(3)}(0) e^{c_{13} t}
+ \int_0^t e^{c_{13}(t-s)} \cdot \bigg\{ c_{10} \io |\na u_{tt}(\cdot,s)|^2 + c_{11} \io |\Del u_t(\cdot,s)|^2 + c_{12}
\bigg\} ds \nn\\
&\le& c_{14} + c_{15} \int_0^t \io |\na u_{tt}|^2 + c_{15} \int_0^t \io |\Del u_t|^2
\qquad \mbox{for a.e.~} t\in (0,\tm),
\eea
where $c_{14}:=\big\{ \sup_{\eps>0} \F_\eps^{(3)}(0) + c_{12} \tm \big\} \cdot e^{c_{13}\tm}$ and
$c_{15}:=(c_{10}+c_{11}) \cdot e^{c_{13}\tm}$ are both finite due to our hypothesis that $\tm<\infty$, and thanks to the fact
that by (\ref{12.22}), the Cauchy--Schwarz inequality, (\ref{99.1}) and (\ref{99.2}),
\bas
\F_\eps^{(3)}(0)
&\le& \frac{\tau}{2} \io \big| \na \big\{ (1+\eps A)^{-1} u_2\big\} \big|^2
+ \Big(\frac{\beta}{2}+1\Big) \io \big| \Del \big\{ (1+\eps A)^{-1}  u_1\big\} \big|^2 \\
& & + \Big(\frac{B}{2} + \frac{\gamma}{4}\Big) \io \big| \Del \big\{ (1+\eps A)^{-1} u_0\big\} \big|^2  \\
&\le& \frac{\tau}{2} \io |\na u_2|^2
+ \Big(\frac{\beta}{{2}}+1\Big) \io |\Del u_1|^2
+ \Big(\frac{B}{2}+\frac{\gamma}{4}\Big) \io |\Del u_0|^2
\qquad \mbox{for all } \eps>0.
\eas
But since (\ref{99.4}) and (\ref{99.3}) in conjunction with (\ref{99.1}), (\ref{99.2}), and (\ref{reg_S}) readily imply that
for a.e.~$t\in (0,\tm)$ we have
\bas
\na \big\{ (1+\eps A)^{-1} u_{tt}\big\} \wto \na u_{tt}
\ \mbox{in } L^2(\Om)
\quad \mbox{and} \quad
\Del \big\{ (1+\eps A)^{-1} u_t\big\} {\wto \Delta u_{t}}
\quad \mbox{in } L^2(\Om)
\eas
as $\eps\searrow 0$, by lower semicontinuity of the norm in $L^2(\Om)$ with respect to weak convergence the inequality in
(\ref{12.15}) together with (\ref{12.3}) entails that the function $y\in L^\infty_{loc}([0,\tm))$ given by
\bas
y(t):=\frac{\tau}{2} \io |\na u_{tt}(\cdot,t)|^2 
+ \frac{\beta}{4} \io |\Del u_t(\cdot,t)|^2,	
\qquad t\in (0,\tm),
\eas
satisfies
\bas
y(t)
&\le& \liminf_{\eps\searrow 0} \F_\eps^{(3)}(t) \\
&\le& c_{14} + c_{15} \int_0^t \io |\na u_{tt}|^2 + c_{15} \int_0^t \io |\Del u_t|^2 \\
&\le& c_{14} + \Big(\frac{2c_{15}}{\tau} + \frac{4c_{15}}{\beta}\Big) \cdot \int_0^t y(s) ds
\qquad \mbox{for a.e.~} t\in (0,\tm).
\eas
Through a Gronwall inequality, again thanks to the finiteness of $\tm$ this implies that actually $y\in L^\infty((0,\tm))$,
and that thus both (\ref{12.1}) and (\ref{12.2}) hold.
\qedd
\mysection{Uniqueness. Proof of Theorem \ref{theo13}} \label{sec:Uniqueness}
For completeness, let us include a brief proof of the uniqueness property claimed in Theorem \ref{theo13}.
\begin{lem}\label{lem14}
	If (\ref{init0}) holds, then for any $T\in (0,\infty]$ there exists at most one strong solution of (\ref{0}) 
	in $\Om\times (0,T)$.
\end{lem}
\proof
Assuming $u$ an $\ou$ to be two strong solutions of (\ref{0}) in $\Om\times (0,T)$, we fix $T_0\in (0,T)$ and then infer from
the local Lipschitz continuity of $f$ and the boundedness of $u$ and $\ou$ in $\Om\times (0,T_0)$ 
that there exists $c_1=c_1(T_0)>0$ such that
\bas
|f(u)-f(\ou)| \le c_1 |u-\ou|
\qquad \mbox{in } \Om\times (0,T_0).
\eas
Thus, if for $(x,t)\in\bom\times [0,T_0]$ we let
\bas
U(x,t):=u(x,t)-\ou(x,t),
\qquad \mbox{and} \qquad
V(x,t):=\int_0^t u(x,s) ds - \int_0^t \ou(x,s) ds
\eas
as well as
\bas
W(x,t):=\int_0^t \int_0^s u(x,\sig)d\sig ds - \int_0^t \int_0^s \ou(x,\sig) d\sig ds,
\eas
then from (\ref{reg_S}) we obtain that, with $B$ taken from (\ref{B}),
\bas
\F(t):=\frac{\tau}{2} \io U^2(\cdot,t)
+ \frac{\beta}{2} \io |\na V(\cdot,t)|^2
+ \gamma \io \na V(\cdot,t) \cdot\na W(\cdot,t)
+ \frac{B}{2} \io |\na W(\cdot,t)|^2,
\qquad t\in [0,T_0],
\eas
defines a function $\F\in C^1([0,T_0])$ which due to (\ref{02}) and Young's inequality satisfies
\bas
\F'(t)
&=& \tau \io UU_t 
+ \beta \io \na U\cdot\na V
+ \gamma \io \na U\cdot\na W
+ \gamma \io |\na V|^2
+ B \io \na V\cdot\na W \\
&=& \io U\cdot \big\{ \tau U_t - \beta\Del V - \gamma \Del W\big\}
+ \gamma \io |\na V|^2
+ B \io \na V\cdot\na W \\ 
&=& \io U \cdot\big\{ -\al U + f(u)-f(\ou)\big\}
+ \gamma \io |\na V|^2
+ B \io \na V\cdot\na W \\
&\le& (|\al|+c_1) \io U^2
+ (\gamma+1) \io |\na V|^2
+ \frac{B^2}{4} \io |\na W|^2
\qquad \mbox{for all } t\in (0,T_0).
\eas
Since, also by Young's inequality,
\bas
\F(t)\ge \frac{\tau}{2} \io U^2
+ \frac{\beta}{4} \io |\na V|^2
+ \frac{B}{4} \io |\na W|^2
\qquad \mbox{for all } t\in (0,T_0),
\eas
this implies that if we let $c_2\equiv c_2(T_0):= \frac{2(|\al|+c_1)}{\tau} + \frac{4(\gamma+1)}{\beta}+B$, then
\bas
\F'(t) \le c_2 \F(t)
\qquad \mbox{for all } t\in (0,T_0).
\eas
The observation that $\F(0)=0$ thus leads to the conclusion that $\F\equiv 0$ and hence $U\equiv 0$ in $\Om\times (0,T_0)$,
which on taking $T_0\nearrow T$ indeed confirms that $u-\ou\equiv 0$ in $\Om\times (0,T)$.
\qedd
Our main result on local existence and uniqueness, and on exclusion of derivative blow-up, has thus been established:\abs
\proofc of Theorem \ref{theo13}. \
It is sufficient to combine Lemma \ref{lem8} with Lemma \ref{lem14} as well as Lemma \ref{lem10}, Lemma \ref{lem11},
and Lemma \ref{lem12}.
\qedd
\mysection{Blow-up. Proof of Theorem \ref{theo45}} \label{Sec:Blow-up}
Our detection of finite-time blow-up in (\ref{0}) will be launched by the following observation on the evolution
of the principal Fourier modes associated with solutions. Again, our analysis in this direction concentrates on the 
integrated problem (\ref{02}).
\begin{lem}\label{lem4}
	Assume (\ref{f}) and (\ref{init}), and let $u$ be taken from Theorem \ref{theo13}. 
	Then, with
	\be{kappa}
	\kappa:=\io e_1,
	\ee
	letting
	\be{y}
	y(t):=\frac{1}{\kappa} \io u(x,t) e_1(x) dx,
	\qquad t\in [0,\tm),
	\ee
	defines a function $y\in C^1([0,\tm))$ which satisfies
	\bea{4.1}
	\tau y'(t) + \al y(t) 
	&=& - \beta \lam_1 \int_0^t y(s) ds
	- \gamma \lam_1 \int_0^t \int_0^s y(\sig) d\sig ds \nn\\
	& & + \frac{1}{\kappa} \io f(u(\cdot,t)) e_1
	+ \frac{t}{\kappa} \io z_1 e_1
	+ \frac{1}{\kappa} \io z_2 e_1
	\quad \mbox{for all } t\in (0,\tm),
	\eea
	where $z_1$ and $z_2$ are as in (\ref{z1}) and (\ref{z2}).
\end{lem}
\proof
The claimed regularity property directly results from (\ref{reg13}), and using (\ref{02}) we compute
\bea{4.2}
\tau y'(t) + \al y(t)
&=& \frac{1}{\kappa} \io \big\{ {\beta \Del v +} \gamma \Del {w}+ f(u) + t z_1 + z_2 \big\} \cdot e_1 \nn\\
&=&  { \frac{\beta}{\kappa} \io \Del v \cdot e_1
	+} \frac{\gamma}{\kappa} \io \Del w \cdot e_1 \nn\\
& & + \frac{1}{\kappa} \io f(u) e_1
+ \frac{t}{\kappa} \io z_1 e_1
+ \frac{1}{\kappa} \io z_2 e_1
\qquad \mbox{for all } t\in (0,\tm).
\eea
As $-\Del e_1=\lam_1 e_1$ in $\Om$ with $e_1|_{\pO}=0$, in line with our definitions of $v$ and $w$ we can here integrate by parts
to rewrite
\bas
\frac{\beta}{\kappa} \io \Del v\cdot e_1
= \frac{\beta}{\kappa} \io v\Del e_1
= - \frac{\beta \lam_1}{\kappa} \io v e_1
= - \frac{\beta \lam_1}{\kappa} \io \bigg\{ \int_0^t u(\cdot,s) ds \bigg\} \cdot e_1
= - \beta\lam_1 \int_0^t y(s) ds
\eas
and, similarly,
\bas
\frac{\gamma}{\kappa} \io \Del w \cdot e_1
= - \frac{\gamma\lam_1}{\kappa} \io \bigg\{ \int_0^t \int_0^s u(\cdot,\sig) d\sig ds \bigg\} \cdot e_1
= - \gamma \lam_1 \int_0^t \int_0^s y(\sig) d\sig ds
\eas
for all $t\in (0,\tm)$, so that (\ref{4.1}) becomes equivalent to (\ref{4.2}).
\qedd
Now, under the additional assumptions on convexity and superlinear growth of $f$ underlying Theorem \ref{theo45}, 
for suitably chosen initial data, the integro-differential equation in (\ref{4.1}) can be turned into a first-order ODI
that must cease to hold by a favorably small time.
\begin{lem}\label{lem44}
	Suppose that (\ref{f}), (\ref{f1}), (\ref{f2}), and (\ref{f3}) hold.
	Then given any $T_0>0$, one can find $K_0=K_0(T_0)>0$ with the property that whenever $u_0\in {W^{4,2}}(\Om)\cap W_0^{1,2}(\Om)$
	is such that (\ref{45.1}) holds, there exists $K_1=K_1(u_0)>0$ such that if
	$u_1\in W^{2,2}(\Om)\cap W_0^{1,2}(\Om)$ satisfies (\ref{45.2}), then it is possible to choose $K_2=K_2(u_1)>0$
	in such a way that for arbitrary $u_2\in W^{2,2}(\Om)\cap W_0^{1,2}(\Om)$ fulfilling (\ref{45.3}),
	for the solution $u$ of (\ref{0}) from Theorem \ref{theo13} we have  
	\be{44.4}
	\tm \le T_0.
	\ee
\end{lem}
\proof
Let $T_0>0$. Then according to (\ref{f3}), we can fix $\xi_1=\xi_1(T_0)>\xi_0$ such that
\be{44.5}
\int_{\xi_1}^\infty \frac{d\xi}{f(\xi)} < \frac{T_0}{4\tau},
\ee
while (\ref{f2}) enables us to find some $\xi_2=\xi_2(T_0)>\xi_0$ fulfilling
\be{44.8}
|\al| \cdot \frac{\xi}{f(\xi)} \le \frac{1}{4}
\qquad \mbox{for all } \xi\ge \xi_2
\ee
and 
\be{44.6}
\beta\lam_1 T_0 \cdot \frac{\xi}{f(\xi)} \le \frac{1}{4}
\qquad \mbox{for all } \xi\ge \xi_2
\ee
as well as
\be{44.7}
\frac{\gamma\lam_1 T_0^2}{2} \cdot \frac{\xi}{f(\xi)} \le \frac{1}{4}
\qquad \mbox{for all } \xi\ge \xi_2.
\ee
We then let
\be{44.9}
K_0 \equiv K_0(T_0) := 
{ \kappa \cdot 
}
\max \{\xi_1 \, , \, \xi_2\},
\ee
and assuming 
{
	$u_0\in W^{4,2}(\Om)\cap W_0^{1,2}(\Om)$
}
to satisfy (\ref{45.1}), we take $K_1=K_1(u_0)>0$ such that
\be{44.10}
\tau K_1 \ge - \al \io u_0 e_1 + \io f(u_0) e_1.
\ee
Given $u_1\in W^{2,2}(\Om)\cap W_0^{1,2}(\Om)$ such that (\ref{45.2}) holds, we finally fix $K_2=K_2(u_1)>0$ in such a way that
\be{44.11}
\tau K_2 \ge - \al \io u_1 e_1
+ \beta \io \Del u_0 \cdot e_1
+ \io f'(u_0) u_1 e_1,
\ee
and henceforth suppose that $u_2\in W^{2,2}(\Om)\cap W_0^{1,2}(\Om)$ satisfies (\ref{45.3}).\abs
To see that then the corresponding solution $u$ has the property in (\ref{44.4}), we assume on the contrary that we had
$\tm>T_0$, and let $y$ be as introduced in Lemma \ref{lem4}, observing that in view of (\ref{45.2}) and (\ref{44.10}),
on the right-hand side of (\ref{4.1}) we can estimate
\bas
\io z_2 e_1
&=& \tau \io u_1 e_1
+ \al \io u_0 e_1
- \io f(u_0) e_1 \\
&\ge& \tau K_1 
+ \al \io u_0 e_1
- \io f(u_0) e_1 \\[2mm]
&>& 0.
\eas
Similarly, (\ref{45.3}) together with (\ref{44.11}) implies that also
\bas
\io z_1 e_1
&=& \tau \io u_2 e_1
+ \al \io u_1 e_1
- \beta \io \Del u_0 \cdot e_1
- \io f'(u_0) u_1 e_1 \\
&\ge& \tau K_2
+ \al \io u_1 e_1
- \beta \io \Del u_0 \cdot e_1
- \io f'(u_0) u_1 e_1 \\[2mm]
&>& 0,
\eas
whence (\ref{4.1}) implies the inequality
\bea{44.111}
\tau y'(t)
> - \al y(t) 
+ \frac{1}{\kappa} \io f(u(\cdot,t)) e_1
- \beta\lam_1 \int_0^t y(s) ds
- \gamma\lam_1 \int_0^t \int_0^s y(\sig) d\sig ds
\eea
for all $t\in (0,\tm)$.
We can now {make use of} the convexity assumption in (\ref{f1}) to see that thanks to the nonnegativity of $e_1$ and the $L^1$
normalization $\io \frac{e_1}{\kappa} =1$ asserted by (\ref{kappa}), the Jensen inequality applies so as to ensure that
\bas
f(y(t))
= f\bigg( \io u(\cdot,t) \cdot \frac{e_1}{\kappa}\bigg)
\le \io f(u(\cdot,t)) \cdot \frac{e_1}{\kappa}
\qquad \mbox{for all } t\in (0,\tm).
\eas
Therefore, (\ref{44.111}) entails that
\be{44.12}
\tau y'(t)
> f(y(t)) - \al y(t)
- \beta\lam_1 \int_0^t y(s) ds
- \gamma\lam_1 \int_0^t \int_0^s y(\sig) d\sig ds
\qquad \mbox{for all } t\in (0,\tm),
\ee
and we claim that our above preparations ensure that, in fact,
\be{44.13}
\tau y'(t) \ge \frac{1}{4} f(y(t))
\qquad \mbox{for all } t\in (0,T_0).
\ee
To verify this, we let
\be{44.14}
\wh{T}_0:=\sup \Big\{ T\in (0,T_0) \ \Big| \ 
{
	y(t)> \frac{K_0}{\kappa}
} 
\mbox{ and } y'(t)> 0 \mbox{ for all } t\in (0,T) \Big\},
\ee
and note that $\wh{T}_0$ is well-defined and positive, because $y$ and $y'$ are continuous on $[0,\tm)$ by Lemma \ref{lem4},
because $y(0)>{\frac{K_0}{\kappa}}$ by (\ref{45.1}), and because, as a particular consequence of (\ref{44.8}),
\bas
\tau y'(0)
> f(y(0)) - \al y(0)
= f(y(0)) \cdot \Big\{ 1 - \frac{\al y(0)}{f(y(0))} \Big\}
\ge f(y(0)) \cdot \frac{3}{4}>0.
\eas
We next take full advantage of (\ref{44.8}) to see that since 
{
	{$y\ge \frac{K_0}{\kappa}\ge \xi_2$} on $(0,\wh{T}_0)$ by (\ref{44.14}) and
	(\ref{44.9}),
	\be{44.15}
	{\frac{\al y(t)}{f(y(t))} 	
		\le |\al| \cdot \frac{y(t)}{f(y(t))}
		\le \frac{1}{4} }
	\qquad \mbox{for all } t\in (0,\wh{T}_0),
	\ee
	whereas additionally using that $y'\ge 0$ on $(0,\wh{T}_0)$ we can rely on (\ref{44.6}) to estimate
	\bea{44.16}
	\frac{\beta\lam_1 \int_0^t y(s) ds}{f(y(t))}
	&\le& \frac{\beta\lam_1 \int_0^t y(t) ds}{f(y(t))} \nn\\
	&=& \beta\lam_1 t \cdot \frac{y(t)}{f(y(t))} \nn\\
	&\le& \beta\lam_1 T_0 \cdot \frac{y(t)}{f(y(t))} \nn\\
	&\le& \frac{1}{4}
	\qquad \mbox{for all } t\in (0,\wh{T}_0),
	\eea
	because $\wh{T}_0 \le T_0$. Likewise,
	{ 
		(\ref{44.7}) ensures that
	}
	\bas
	\frac{\gamma\lam_1 \int_0^t \int_0^s y(\sig) d\sig ds}{f(y(t))}
	&\le& \frac{\gamma\lam_1 \int_0^t \int_0^s y({t}) d\sig ds}{f(y(t))} \\
	&=& \frac{\gamma\lam_1 {t^2}}{2} \cdot \frac{y(t)}{f(y(t))} \\
	&\le& \frac{\gamma\lam_1 {T^2_0}}{2} \cdot \frac{y(t)}{f(y(t))} \\
	&\le& \frac{1}{4}
	\qquad \mbox{for all } t\in (0,\wh{T}_0),
	\eas
	which along with (\ref{44.15}) and (\ref{44.16}) shows that (\ref{44.12}) implies the inequality
	\be{44.17}
	\tau y'(t) \ge \frac{1}{4} f(y(t))
	\qquad \mbox{for all } t\in (0,\wh{T}_0).
	\ee
	As $f$ is positive on 
	{
		$[\frac{K_0}{\kappa},\infty)$,
	}
	this especially ensures that actually $\wh{T}_0=T_0$, for otherwise the continuity of $y'$
	would entail that, on the one hand, we should either have 
	{
		$y(\wh{T}_0)=\frac{K_0}{\kappa}$ 
	}
	or $y'(\wh{T}_0)=0$, and that, on the other hand,
	$y'>0$ on $[0,\wh{T}_0]$ and hence both 
	{
		$y(\wh{T}_0)>\frac{K_0}{\kappa}$
	} 
	and $y'(\wh{T}_0)>0$.
	Since this is impossible, we indeed obtain (\ref{44.13}) from (\ref{44.17}), and thereupon the claim simply results
	by integration:
	From (\ref{44.13}) and the inequality 
	{
		$y(0)\ge \frac{K_0}{\kappa}$,
	}
	namely, we obtain that
	\bas
	\frac{t}{4\tau}
	\le \int_{y(0)}^{y(t)} \frac{d\xi}{f(\xi)}
	\le 
	{
		\int_{\frac{K_0}{\kappa}}^\infty \frac{d\xi}{f(\xi)}
	}
	\qquad \mbox{for all } t\in [0,T_0],
	\eas
	so that by (\ref{44.9}) and (\ref{44.5}), 
	\bas
	\frac{t}{4\tau} \le \int_{\xi_1}^\infty \frac{d\xi}{f(\xi)} < \frac{T_0}{4\tau}
	\qquad \mbox{for all } t\in [0,T_0],
	\eas
	which when evaluated at $t=T_0$ leads to a contradiction.
	Our hypothesis that $\tm$ be larger than $T_0$ must therefore have been false.
	\qedd
	In conclusion, we have thereby achieved our main result on finite-time $L^\infty$ blow-up in (\ref{0}):\abs
	\proofc of Theorem \ref{theo45}. \quad
	The claim immediately results if we combine Lemma \ref{lem44} with Theorem \ref{theo13}.
	\qedd

	\bigskip

	{\bf Acknowledgments.} \
	The second author acknowledges support of the {\em Deutsche Forschungsgemeinschaft} (Project No.~444955436).
	The authors declare that they have no conflict of interest.

	\bigskip
	
	{\bf Data availability statement.} \quad
	Data sharing is not applicable to this article as no datasets were generated or analyzed during the current study.

\bibliography{references}
\bibliographystyle{siam} 
\end{document}